\documentclass[leqno,11pt]{article}%

\usepackage{stmaryrd}
\usepackage{setspace,bbm,mathrsfs}

\usepackage[pagebackref,hypertexnames=false]{hyperref}
\usepackage{bbm,subfigure}
\usepackage{verbatim,amssymb,amsfonts,tikz-cd}
\usepackage{multicol,xfrac}
\usepackage{fullpage,tikz}

 \usepackage{amsmath,amsthm,bbm,amssymb}

\DeclareMathOperator{\supp}{supp}
\DeclareMathOperator{\im}{im}

\DeclareMathOperator{\vol}{Vol}

\DeclareMathOperator{\PH}{PH}

\def\R{\mathbb{R}}

\def\vsimp{V_{\mathrm{simp}}}
\def\Psir{\Psi_{\mathrm{rad}}}
\def\Psis{\Psi_{\mathrm{sph}}}

\def\S{S}

\def\cC{\mathcal{C}}
\def\cD{\mathcal{D}}
\def\cF{\mathcal{F}}

\def\cH{\mathcal{H}}
\def\cI{\mathcal{I}}

\def\cL{\mathcal{L}}

\def\cP{\mathcal{P}}

\def\cR{\mathcal{R}}

\def\cS{U}

\def\cX{\mathcal{X}}
\def\cY{\mathcal{Y}}

\def\cL{\mathcal{L}}

\newcommand{\E}{\mathbb{E}} 

\newcommand{\given}{\;|\;}
\newcommand{\mean}[1] {\E\left\{{#1}\right\}}

\newcommand{\meanx}[1] {\E\{{#1}\}}
\newcommand{\cmean}[2] {\E\left\{#1\given #2\right\}}

\newcommand{\cmeanx}[2] {\E\{#1\given #2\}}
\newcommand{\ind}{\boldsymbol{\mathbbm{1}}} 
\newcommand{\indf}[1]{\ind\set{#1}} 
\newcommand{\indfx}[1]{\ind\{{#1}\}} 

\newcommand{\iprod}[1]{\left\langle{#1}\right\rangle}
\newcommand{\set}[1]{\left\{#1\right\}}

\newcommand{\param}[1]{\left(#1\right)}
\newcommand{\abs}[1] {\left| {#1}\right|}
\newcommand{\floor}[1] {\left\lfloor{#1}\right\rfloor}
\newcommand{\ceil}[1] {\left\lceil{#1}\right\rceil}
\newcommand{\prob}[1]{\mathbb{P}\left(#1\right)}
\newcommand{\probx}[1]{\mathbb{P}(#1)}
\newcommand{\cprob}[2]{\mathbb{P}\left(#1\given #2\right)} 
\newcommand{\cprobx}[2]{\mathbb{P}(#1\given #2)} 

\newcommand{\eps}{\varepsilon}
\newcommand{\bz}{\mathbf{z}}
\newcommand{\by}{\mathbf{y}}
\newcommand{\bx}{\mathbf{x}}

\newcommand{\rmax}{r_{\max}}
\newcommand{\lmax}{\Lambda_{\max}}

\def\bth{\boldsymbol\theta}
\providecommand{\setthms}[1]{#1}
\setthms{
\newtheorem{lem}{Lemma}[section]
\newtheorem{thm}[lem]{Theorem}
\newtheorem{prop}[lem]{Proposition}
\newtheorem{cor}[lem]{Corollary}

\newtheorem{rem}[lem]{Remark}

\theoremstyle{definition}
\newtheorem{defn}[lem]{Definition}

}

\def\Hg{\mathrm{H}}
\def\PHg{\mathrm{PH}}

\def\lk{\mathrm{lk}}

\newcommand{\ninf}{n\to\infty}
\newcommand{\pois}[1]{\mathrm{Poisson}\param{{#1}}}

\newcommand{\fmax}{f_{\max}}
\newcommand{\fmin}{f_{\min}}

\newcommand{\limninf}{\lim_{\ninf}}

\newcommand{\bs}{\backslash}
\definecolor{mygreen}{rgb}{0, 0.68, 0.31}
\definecolor{myred}{rgb}{1.0, 0,0}

\numberwithin{equation}{section}

\def\bsplit#1\esplit{\begin{split} #1 \end{split} }
\def\splitb#1\splite{\begin{split} #1 \end{split} }
\def\beq#1\eeq{\begin{equation} #1 \end{equation}}
\def\eqb#1\eqe{\begin{equation} #1 \end{equation}}

\def\dgm{\mathrm{dgm}}

\newcommand{\cK}{\mathcal{K}}
\newcommand{\binmp}[2] {\mathrm{Binomial}({#1};{#2})}
\newcommand{\poisp}[2] {\mathrm{Poisson}({#1};{#2})}
\newcommand{\pairing}{\phi}
\newcommand{\Cl}{\mathrm{cl}}
\newcommand{\clStr}{\mathrm{st}}
\newcommand{\Csc}{\mathrm{Csc}}
\newcommand{\VR}{\cR}
\newcommand{\mcech}{\cC}
\def\pointset{\cX}

\newcommand{\birth}{\mathrm{birth}}
\newcommand{\death}{\mathrm{death}}
\def\nn{\nu}
\newcommand{\birthx}{{b}}
\newcommand{\deathy}{{d}}

\DeclareMathOperator{\coker}{coker}
\DeclareMathOperator{\rk}{rk}

\newcommand{\cone}{\mathrm{Cone}}
\newcommand{\origin}{\mathbf{0}}
\title{Universality in Random Persistent Homology and Scale-Invariant Functionals}
\author{Omer Bobrowski\footnote{o.bobrowski@qmul.ac.uk} \ and  Primoz Skraba\footnote{p.skraba@qmul.ac.uk}}
\date{}

\begin{document}

\maketitle

\begin{abstract}
In this paper, we prove a universality result for the limiting distribution of persistence diagrams arising from geometric filtrations over random point processes. Specifically, we consider the distribution of the ratio of persistence values  (death/birth), and  show that for fixed  dimension, homological degree 
 and filtration type (\v Cech or Vietoris-Rips), the limiting distribution is independent of the underlying point process distribution, i.e., universal.  
 In proving this result, we present a novel general framework for  universality in scale-invariant functionals on point processes. Finally, we also provide a number of new results related to Morse theory in random geometric complexes, which may be of an independent interest. 
\end{abstract}

\section{Introduction}
Persistent homology has established itself  in the last two decades as an important tool both within theoretical mathematics 
\cite{alpert_configuration_2024,ellis_persistent_2011,entov_legendrian_2022,
polterovich_topological_2020} and in the field of topological data analysis (TDA) \cite{carlsson_topological_2021,chazal_high-dimensional_2017,dey_computational_2022,edelsbrunner_computational_2010}. In its most prevalent form, persistent homology is an invariant associated with a filtration of topological spaces $\cF = \set{X_t}_{t\in \R}$. Applying the homology functor, we obtain a \emph{persistence module} -- a sequence of homology groups associated to each $t\in \R$ along with homomorphisms $\psi_{s,t}: \Hg_* (X_s)\rightarrow \Hg_*(X_t)$ for all $s\leq t$
.
Taking homology over a field,  the persistent homology module admits a simple decomposition, formalizing the notion of \emph{tracking changes} in homology over the filtration, with new generators appearing (birth) and later becoming trivial (death). The collection of all (birth,death) pairs for classes appearing throughout the filtration, can be considered as a point-set, or a discrete measure, in $\R^2$, known as the \emph{persistence diagram}.

In this paper we focus on filtrations of  simplicial complexes generated by random point-processes in $\R^d$ (indexed by a radius parameter $r$).
Extensive study over the past two decades has revealed  novel phase transitions describing major changes related to persistent homology \cite{bobrowski_homological_2022-1,bobrowski_homological_2022,muller_contractibility_2022}, as well as various limit theorems \cite{divol_density_2019,hiraoka_limit_2018,krebs_asymptotic_2019,owada_convergence_2020,yogeshwaran_topology_2015}. A key property of all previous results, is that the limiting behavior depends on the underlying distribution of the point-process.

Viewing persistence diagrams as Radon measures in $\R^2$ (the birth-death distribution), it was shown in \cite{divol_density_2019,hiraoka_limit_2018} that for random geometric complexes, these measures have a deterministic limit. The limiting measure, however, depends on the underlying point-process distribution in $\R^d$.
Contrary to this, 
\cite{bobrowski_universal_2023} presented comprehensive experimental evidence leading to a number of surprising conjectures.
The main claim is that the distribution of  the ratio values $\pi = \death/\birth$ is \emph{universal} -- independent of the underlying point-process distribution. 

\begin{figure}
    \centering
\includegraphics[width=0.8\textwidth]{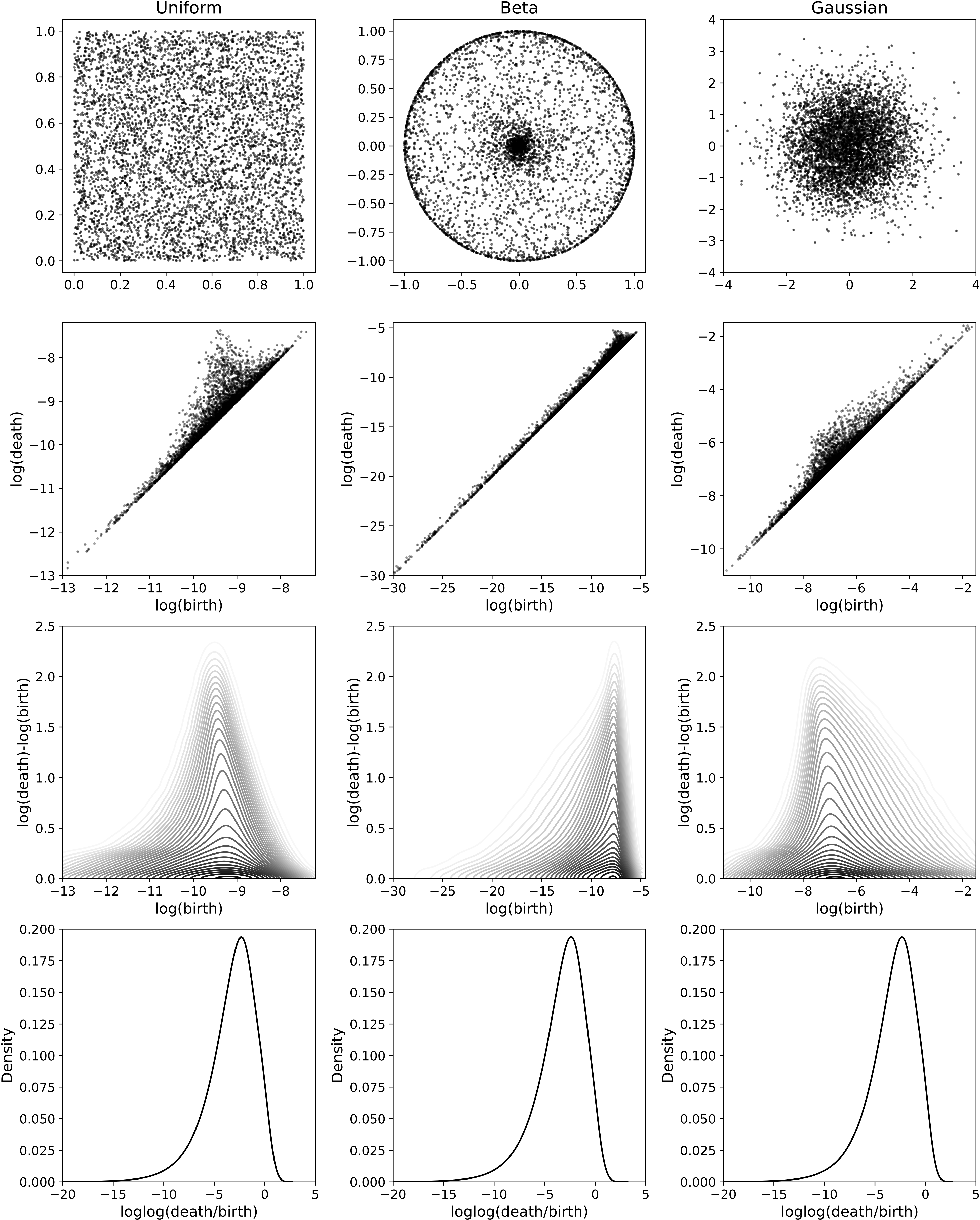}
    \caption{Universal distribution for the $\pi$-values. 
    We use iid samples of size $n=5000$, the \v Cech filtration, and homological degree $k=1$. Each column examines a different distribution.
    {\bf Left:} Uniform distribution in $[0,1]^2$; {\bf Center:} Uniform distribution for the angle and beta distribution for the radius; {\bf Right:} Standard normal distribution.
    The first and second rows present the original point sets and their corresponding persistence diagrams 
    (shown in log-scale). The third row presents the contour lines for the estimated two-dimensional density of the diagram (rotated by 45 degrees).
    The last row presents the estimated density for the $\pi$-values (loglog scale). 
    While the resulting distributions in the second and third rows are vastly different between the three point sets, the distribution of the $\pi$-values is identical.}
    \label{fig:example}
\end{figure}

In this paper, we prove the first conjecture in \cite{bobrowski_universal_2023} -- if we fix the dimension of the space $d$, the homological degree $k$, and the construction type (i.e., \v Cech or Vietoris-Rips), the empirical distribution of the $\pi$-values, generated by iid point-processes is universal --  see Figure~\ref{fig:example} for an example. 
This  conjecture was named
``Weak Universality'', 
since the experiments in \cite{bobrowski_universal_2023} suggest a much stronger universal behavior (i.e., beyond iid,  independent of dimensions, homological degree, and complex-type). Aside from it being a surprising result in the context of stochastic topology, universality in persistence diagram has immense potential in developing statistical applications in TDA. For example, it can lead to robust topological inference methods \cite{bobrowski_universal_2023}, new and powerful clustering algorithms \cite{bobrowski_cluster_2023}, dimensionality estimation, measuring system disorder, and more. In particular, it implies that statistical-topological tests can be designed with very little knowledge about the nature of the data.

To prove universality in persistence diagrams, 
we  develop a novel generic framework for a \emph{universal law of large numbers} that applies to scale-invariant point-process functionals. 
Let $\cX_n = \set{X_1,\ldots,X_n}$ be a collection of $n$ iid random variables, with a probability density function $f:\R^d\to \R$. We show that under some conditions, for scale-invariant and linearly scaling functionals $\cH$, we have $\frac{1}{n} \cH(\cX_n) \xrightarrow{a.s.} \cH^*$ where $\cH^*$ is a (nonzero) deterministic constant independent of $f$ (i.e., universal). 
While our main interest is for persistence diagrams, 
we also 
demonstrate the broader potential of this framework using the example of the  degree distribution in random k-nearest neighbor graphs. 

\vspace{5pt}
\noindent{\bf Related work.}
Since the early days of persistent homology, there has been an ongoing effort to reveal the probabilistic behavior of random persistence diagrams. The seminal result in \cite{hiraoka_limit_2018} shows that for stationary point processes, random persistence diagrams have a deterministic limiting measure (in the weak sense). This result was further explored in numerous directions, including the existence of density \cite{divol_density_2019}, central limit theorems \cite{krebs_asymptotic_2019}, special limits for the sparse regime \cite{owada_convergence_2022}, extreme value analysis \cite{owada_convergence_2020}, spanning acycles \cite{skraba_central_2023}, marked point processes \cite{shirai_limit_2022}, generic weight functions \cite{divol_choice_2019}, and dependencies \cite{krebs_limit_2021}. In  these results, and more generally in stochastic geometry and topology, the limiting distributions and related constants depend on the underlying point-process distribution. Our results show that by projecting the persistence measures from $\R^2$ to $\R$, via the $(\death/{\birth})$ transformation, we obtain a random measure whose (non-random) limit is independent of the underlying distribution, and hence universal.

The first (and to the best of our knowledge, the only) description of universal limits for point-process functionals appeared in \cite{penrose_weak_2003}. The goal there is to prove weak laws of large numbers for functionals over binomial processes, that can be written as sums of score functions for individual points. The limit provided in \cite{penrose_weak_2003} has an explicit formula in terms of the density $f$, and in the case where the functional is scale-invariant, the dependency in $f$ vanishes. The examples presented there are related to the k-th nearest neighbor (k-NN) graph and Voronoi/Delaunay graphs. The functionals considered include the number of components, sub-graph counts, and sums of edge-lengths. While that framework is quite powerful, we are not able to apply it for persistence diagrams due to three main challenges: (a) Persistent-cycles counting is not naturally expressed as a sum of score functions; (b) The coupling used in \cite{penrose_weak_2003}, presents a limitation of the results to the thermodynamic limit, which in our case is essential to bypass; (c) For the most part, a required assumption in \cite{penrose_weak_2003} is that the density has a compact and convex support, and is lower-bounded away from zero. This poses a strong limitation which we wish to avoid. The framework we develop here is quite different than the coupling method in \cite{penrose_weak_2003}, and  provides a ``bottom-up'' constructive approach, that reveals the sources of the universal behavior. In addition to addressing (a)-(c), our framework allows for functionals that are scale-invariant only in an asymptotic sense, and provides conditions for strong laws of large numbers.

\section{Preliminaries}

 We begin with a brief review some basic definitions required to state  the main results.

\subsection{Geometric complexes}

Our primary concern will be finite point sets $\pointset\subset \R^d$. Commonly, if these points lie on or near some underlying subspace, we can study the topology of subspace by using discrete approximations based on $\pointset$. In this paper, we consider two of the most common geometric constructions used in this context (see \cite{edelsbrunner_computational_2010}).

\noindent {\bf The \v Cech complex}  at radius $r$ over $\pointset$  is an abstract simplicial complex defined as,
$$\mcech_r(\pointset) = \left\{\mathcal{I} \subset \pointset : \bigcap\limits_{x\in\mathcal{I}  } B_r(x) \neq \emptyset \right\},$$
where $B_r(x)$ is a closed ball of radius $r$, centered at $x$.

\noindent{\bf The Vietoris-Rips complex} at radius $r$  over $\pointset$ is defined as,
$$\VR_r(\pointset) = \left\{\mathcal{I} \subset \pointset P: |x_1-x_2|\le r, \ \forall x_1,x_2\in\mathcal{I}  \right\},$$
where by $|\cdot|$ we refer to the Euclidean norm.

\begin{rem}We use the closed edge definition  (i.e., $|x_1-x_2|\le r$ compared to $|x_1-x_2| < r$) to be consistent with the \v Cech complex and its corresponding Morse theory later. This may result in pathologies, even for compact metric spaces. However,  for the random point processes we study, these do not occur, as for any fixed $r$, the probability of having any pair with $|x_1-x_2|=r$ is zero. 
\end{rem}

\subsection{Persistent homology}

In this paper, we  restrict ourselves to homology with field coefficients.
In this case, given a filtration of spaces $\mathcal{F} = \{X_t\}_{t\in\R}$ (i.e.,  $X_{t_1}\subset X_{t_2}$ for $t_1\le t_2$), the  persistent homology $\PHg_k(\cF)$  is a module of vector spaces  
$\{\Hg_k(X_t)\}_{t\in \R}$ along with linear maps  $\Hg_k(X_{t_1}) \rightarrow \Hg_k(X_{t_2})$ for all $t_1\le t_2$ induced by the corresponding inclusion maps. 
If $\Hg_k(X_t)$ is finite dimensional for every $t$, then  the module $\PH_k(\cF)$ admits a simple decomposition into interval summands (i.e., rank-one submodules which are non-trivial over an interval of parameters $[b,d)$, and zero everywhere else, denoted $\mathbb{I}[b,d)$). 
\begin{thm}[\cite{crawley-boevey_decomposition_2015}, Theorem 1.1]
Any pointwise finite-dimensional persistence module is
a direct sum of interval summands, up to an isomorphism. In other words,
\eqb\label{eqn:ph_decomp}
\PHg_k(\mathcal{F}) \cong \bigoplus\limits_{i} \mathbb{I}[\birthx_i,\deathy_i).
\eqe
where $\birthx_i$ and $\deathy_i$ are the lower and upper bounds of the interval where the $i$-th summand is non-trivial. 
\end{thm}

A common way to summarize the information encoded in persistent homology, is to look at the collection of $(\birthx_i,\deathy_i)$ pairs in the decomposition. We note that, as is standard in persistent homology,  we only consider  pairs with $b_i<d_i$. 

\begin{defn}
The \emph{persistence diagram} associated with a filtration $\cF$, is a multiset of points in $\R^* \times \R^*$, where $\R^*:= \R\cup\{\infty\}$, defined as 
$$\dgm_k(\mathcal{F}) := \left\{(\birthx_i,\deathy_i),  \;\;i\in \cI :  \PHg_k(\mathcal{F})\cong\bigoplus\limits_{i\in \cI} \mathbb{I}[\birthx_i,\deathy_i)   \right\} $$
\end{defn}
%
We will  restrict our discussion to filtrations generated by the geometric constructions $\cC_r$ and $\cR_r$ defined above, indexed by the distance parameter $r$.
These are often referred to as  \emph{geometric} or \emph{proximity} filtrations.
If $\pointset$ is finite then the homology is finite-dimensional for all values of $r$, and the persistence diagram is well defined. 

A key motivation for using persistent homology in the context of data analysis, is that  $k$-cycles with a long lifetime 
are assumed to represent important topological features. In this paper, we focus on the lifetime of cycle in a multiplicative perspective.
Let $\cF$ be a filtration, and let $p_i=(\birthx_i,\deathy_i)\in \dgm_k(\cF)$ for some $k>0$. We define the corresponding `$\pi$-value', as\
\[
\pi(p_i) := \frac{\deathy_i}{\birthx_i}.
\]
The motivation to study the ratio rather the difference ($\deathy_i-\birthx_i$) to measure cycle-size stems from  the facts that: (a) it is a scale-invariant measure, and (b) it is robust to statistical outliers (see more information in \cite{bobrowski_universal_2023}).
Note that in  geometric filtrations all cycles in degree $0$ are  born at $r=0$. Therefore, we only consider the $k$-th homology for $k>0$. 
\begin{rem}
Throughout this paper, $\Hg_k(\cdot)$ always refers to the reduced homology. As we only consider $k>0$, this does not change the persistence diagrams or functionals we study, but does simplify the exposition and proofs. 
\end{rem}

\subsection{Probabilistic settings}\label{sec:prob_set}

\subsubsection*{Point processes}
Our main statements will be proved for two tightly-related point processes.
The \emph{binomial process}, denoted $\cX_n\sim\binmp{n}{f}$, is generated by taking $n$ iid points with probability density $f$. The \emph{Poisson process} $\cP_n\sim\poisp{n}{f}$ is generated the same way, but  the number of points is $N\sim \pois{n}$ (i.e., random). To avoid confusion, and also to highlight the fact that the rate of the Poisson process does not need to be an integer, we will use $\cP_\nn$ instead of $\cP_n$ throughout the paper.

\subsubsection*{Good densities}

Our goal is to prove the main results for the largest possible  class of distributions. We assume that all such distributions are equipped with a probability density function $f:\R^d\to \R$.
We denote the support of $f$ by $\cS = \mathrm{cl}(f^{-1}(0,\infty))$ (where $\mathrm{cl}$ denotes the closure). 
As the framework requires a number of conditions, our statements for persistent homology will be proved for three types of distributions: 
\begin{itemize}
    \item {\bf Type I:} $\cS$ is compact, $\inf_{\cS} f >0$.
    \item {\bf Type II:} $\cS$ is compact, $\inf_{\cS} f = 0$.\\
    Here, we need an additional assumption. Let $\cS_0 = f^{-1}(0)\cap \cS$. For every $x\in \cS$ define $\delta(x) = \inf_{z\in \cS_0} |x-z|$. We assume that there exist $q,\delta_0,C_1,C_2>0$ such that if $\delta(x)\le \delta_0$ then
    \[
        C_1(\delta(x))^q \le f(x) \le C_2(\delta(x))^q.
    \]
    In other words, the growth rate of the density around its zeros is polynomial.
    \item {\bf Type III:} $\cS = \R^d$.\\
    Here, we will focus on densities of the following type,
    \[
        f(x)= c e^{-\Psi(x)},
    \]
    where $c>0$ is a normalization   constant, and $\Psi$  decomposes into radial and spherical components,
    \[
        \Psi(x) = \Psir(|x|)\Psis(x/|x|).
    \]
    The spherical component is assumed to be a Lipschitz function $\Psis:\S^{d-1}\to [1,\infty)$. The radial part is assumed to be $\Psir:[0,\infty)\to [0,\infty)$, with $\lim_{\tau\to\infty}\Psir(\tau) = \infty$.  We further assume that there exists $R_0>0$ such that for all $\tau\ge R_0$, the following holds:
    \begin{itemize}
        \item[(III.1)] $\Psir(\tau) \ge A_1\log(\tau)$,  for some $A_1>d$.
        \item[(III.2)] $\Psir'(\tau) \ge A_2/\tau$,  for some $A_2>0$.
        \item[(III.3)] $|\Psir(\tau+\eps)-\Psir(\tau)| \le A_3 \eps \Psir(\tau)$ for any small enough $\eps>0$.
    \end{itemize}
\end{itemize}
Note that (III.2) implies that $\Psir$ is increasing in $\tau$ (for $\tau\ge R_0)$.
Type III is seemingly restrictive,  however, it contains many of the well-known distributions. For example,
\begin{itemize}
    \item The normal distribution (with any non-singular covariance matrix).
    \item Any density of the type $f(x)\propto e^{-|x|^\alpha}$ for $\alpha>0$.
    \item Any density of the type $f(x)\propto e^{-e^{|x|^\alpha}}$ for $\alpha\in (0,1)$.
    \item Any density of the type $f(x) = (1+|x|^\alpha)^{-1}$, for $\alpha>d$. 
\end{itemize}

\begin{rem}
In some cases, assumption (III.3) may still be too restrictive (for example $f(x) = e^{-e^{|x|^\alpha}}$ with $\alpha \ge 1$). In this case we can replace it with the following weaker assumption.

\noindent
\emph{(III.3') Let $R = \Psir^{-1}\param{\frac{k+1}{1-d/A_1}\log\nn}$. Then there exists $r=r(\nn)\to 0$ such that $\nn r^d\to\infty$, and for all $\tau \in [R_0,R]$ we have
        \[
        |\Psir(\tau+r)-\Psir(\tau)| \le \frac{A_3}{\log \nn} \Psir(\tau).
        \]}
\end{rem}

\section{Main results}

Let $\pointset$ be a finite subset of $\R^d$ and $\cK = \set{K_r}_{r=0}^\infty$ denote either the \v Cech or Vietoris-Rips filtrations generated by $\pointset$ with  $\dgm_k(\cK)$ the persistence diagram for the $k$-th degree homology.
We define the $k$-th \emph{persistence measure} $\Pi_k$ as follows. For every Borel set $B\subset \R$,
\[
\Pi_k(B) := \abs{\set{p\in\dgm_k(\cK) : \pi(p)\in B}}.
\]
Note that $\Pi_k((-\infty,1]) = 0$, so the actual support of the measure is $[1,\infty)$.
Since $\Pi_k$ is a finite measure, we can turn it into a probability measure $\Psi_k$ by defining
\[
\Psi_k(B) := \frac{\Pi_k(B)}{\Pi_k(\R)}.
\]
In the case where $\cX$ is random, we can study the expected measure $\bar\Pi_k$
\[
\bar\Pi_k(B) := \mean{\Pi_k(B)},
\]
and the associated probability measure $\bar\Psi_k$,
\[\bar\Psi_k(B) := \frac{\bar\Pi_k(B)}{\bar\Pi_k(\R)}.
\]

The main theorems in this paper refer to the persistence diagrams generated by either the Poisson or the binomial processes, where the number of points ($n$) goes to infinity. In this context, it is shown in \cite{kahle_random_2011} that the so-called \emph{thermodynamic limit} (i.e., $nr^d = \mathrm{const}$) is the dominant range where most $k$-cycles are born and die.
Thus, our statements in this paper are focused on covering this entire range. In other words, we will re-define the persistence measure to be 
\eqb\label{eqn:Pi_k}
    \Pi_{k,n}(B) := \abs{\set{p\in\dgm_k(\cK) : \pi(p)\in B, \death(p) \le \rmax}},
\eqe
where $\rmax =\rmax(n)$ is chosen such that $\rmax\to 0$ (as $n\to\infty$) and $n\rmax^d\to\infty$, so that we focus our attention on small cycles, but make sure to include all of those cycles generated within the thermodynamic limit. Considering the dichotomy between topological \emph{signal} and \emph{noise} (cf., \cite{bobrowski_universal_2023}), we can say that our analysis here applies to the noisy cycles (which make the vast majority of cycles in these settings).
The bound on the death time ($\rmax$) will be implicit in the main statements, but will be explicitly used in the proof. Our limiting results are in the sense of the weak convergence of measures \cite{billingsley_convergence_2013}.

\begin{thm}\label{thm:main_poisson}
Let $f:\R^d\to\R$ be a good density function, $\cP_\nn\sim\poisp{\nn}{f}$, and let $\Pi_{k,\nn}, \Psi_{k,\nn}$ be the corresponding  persistence measures. Then,
\[
\lim_{\nn\to\infty} \frac1\nn \bar\Pi_{k,\nn} = \Pi_k^*,
\]
where $\Pi_k^*$ is a finite measure on $\R$. The measure $\Pi_k^*$ depends on $d$ and $k$, as well as the choice of filtration (\v Cech or Vietoris-Rips), but is otherwise independent of $f$. Furthermore, defining the probability measure $\Psi_k^*:={\Pi_k^*}/{\Pi_k^*(\R)}$, we have
\[
\lim_{\nn\to\infty} \bar\Psi_{k,\nn} = \Psi_k^* .
\]
\end{thm}

For the binomial process, we also obtain a strong law of large numbers. 

\begin{thm}\label{thm:main_binom}
Let $f$ be a good density function,  $\cX_n\sim\binmp{n}{f}$, and let $\Pi_{k,n}, \Psi_{k,n}$ be the corresponding  persistence measures. Then almost surely,
\[
\lim_{n\to\infty} \frac1n \Pi_{k,n} = \lim_{n\to\infty} \frac1n \bar\Pi_{k,n} = \Pi_k^*,
\]
and
\[
\lim_{n\to\infty} \Psi_{k,n} = \lim_{n\to\infty} \bar\Psi_{k,n} = \Psi_k^* ,
\]
where $\Pi_k^*,\Psi_k^*$ are the same universal measures as in Theorem \ref{thm:main_poisson}. 
\end{thm}

The remaining sections of this paper are dedicated to proving Theorems \ref{thm:main_poisson} and \ref{thm:main_binom}.
The key property underlying the universality behavior, is the \emph{scale invariance} of the $\pi$-values. By that, we mean that if $\Pi_k$ is the persistence measure for a set $\pointset\subset\R^d$, and $\Pi'_k$ is the persistence measure for $\pointset' = c \pointset$ for some $c>0$, then $\Pi'_k\equiv \Pi_k$.  We therefore begin  with short detour, proving universality in the more general setting of scale-invariant functionals.

\section{Universal limits for scale-invariant functionals}\label{sec:affine_inv}

We first develop a generic framework that can be applied to a wide class of scale-invariant functionals. 
While we originally developed this framework to address the persistence measure, we observed that the conditions required for universality  are not unique to the persistence measure functional.
To demonstrate that, in Section \ref{sec:knn} we prove universal limits for the vertex degree in the random $k$-NN graph. We note that while this example, along with other functionals of the  $k$-NN graph and the Voronoi/Delaunay graphs, were studied in \cite{penrose_weak_2003}, the framework provided here covers functionals that are not expressed as score functions, which is crucial in the case of the persistence measure. Additionally, our framework poses fewer restrictions on the density function, and provides a strong law of large numbers.

\subsection{Universality statements}

Our framework concerns functionals $\cH_\eta$ applied to finite subsets of $\R^d$, that are optionally  controlled by a parameter $\eta$ (e.g., determining the value $\rmax(\eta)$ for the persistence measures). For every such functional, we define  $\cD(\cH)$ as the collection of all probability density functions $f:\R^d\to\R$ for which we want to prove universality. We start by listing the conditions required for our generic universality statements. 

\noindent{\bf
Translation invariance.} Let $\pointset$ be a finite subset of $\R^d$. We say that  $\cH_\eta$ is  \emph{translation invariant} if for every $x_0\in \R^d$ 
\eqb\label{eqn:assum_translate}
    \cH_\eta(\pointset+x_0) = \cH_\eta(\pointset),
\eqe
where $\pointset+x_0 = \set{x+x_0 : x\in \pointset}$.

\noindent{\bf
Linear scale.}
Let $\cP^*_\nn\sim\poisp{\nn}{\ind_{Q^d}}$  with $Q^d = [0,1]^d$, which we name as the \emph{reference model}. Then there exists $\cH^*>0$ such that
\eqb\label{eqn:assum_linear}
    \lim_{\nn\to\infty} \frac1\nn\mean{\cH_\nn(\cP^*_\nn)} = \cH^*.
\eqe

\noindent{\bf
Scale invariance}
Let $\eps>0$, and let $\cP_\nn^{(\eps)} \sim\poisp{\nn}{\eps^{-d}\ind_{\eps Q^d}}$, be a homogeneous Poisson process on $\eps Q^d$ with rate $\nn$. Then
\eqb\label{eqn:assum_scale}
    \lim_{\nn\to\infty} \frac1\nn\meanx{\cH_\nn(\cP^{(\eps)}_\nn)} = \cH^*.
\eqe
In other words, changing the scale of the box does not affect the limit. Note, that if $\cH_\nn$ is strictly scale-invariant (i.e., $\cH_\nn(c\cX) = \cH_\nn(\cX)$ for all $\cX$), then \eqref{eqn:assum_scale} holds trivially.

\noindent{\bf
Homogeneity.}
Let $f\in \cD(\cH)$, and $\cP_\nn\sim\poisp{\nn}{f}$. For any constant $c>1$, 
\eqb\label{eqn:assum_homog}
\lim_{\nn\to\infty} \frac{\mean{\cH_{c\nn}(\cP_\nn)}}{\mean{\cH_\nn(\cP_\nn)}} = 1.
\eqe

\noindent{\bf
Additivity.}
Take any integer $m>1$, and divide $Q^d$ into $M=m^d$ sub-cubes $Q_1,\ldots, Q_M$ of side-length $1/m$. For $i=1,\ldots,M$ let $\cP_{\nn,i}\sim\poisp{c_i\nn}{M\cdot\ind_{Q_i}}$ be independent homogeneous Poisson processes, for some $c_i\ge 0$ with $\sum_ic_i=1$. Define $\cP_\nn = \bigcup_i \cP_{\nn,i}$, then
\eqb\label{eqn:assum_additive}
\lim_{\nn\to\infty} \frac1\nn \abs{\mean{\cH_\nn(\cP_\nn)} - \sum_i \mean{\cH_\nn(\cP_{\nn,i})}} = 0.
\eqe

\noindent{\bf
Continuity.}
Let $f\in \cD(\cH)$, and 
 $\cP_\nn\sim \poisp{\nn}{f}$. 
 Let $f_1,f_2,\ldots$ be a sequence of  functions, such that for all $x\in \R^d$,
\[
\splitb
    0\le f_i(x) &\le f_j(x),\quad \forall i<j,\\
    \lim_{i\to\infty} f_i(x) &= f(x).
\splite
\]
Setting
\[
\delta_i = 1-\int_{\R^d} f_i(x)dx,
\]
we further assume that $\tilde f_i = f_i/(1-\delta_i) \in \cD(\cH)$.
Let $\cP_{\nn,i}\sim\poisp{(1-\delta_i)\nn}{\tilde f_i}$ (i.e., $\cP_\nn$ has intensity function $\nn f_i$). Then,
\eqb\label{eqn:assum_stable}
\lim_{i\to\infty}\lim_{\nn\to\infty} \frac1\nn \abs{\mean{\cH_\nn(\cP_\nn)} - \mean{\cH_\nn(\cP_{\nn,i})}} = 0.
\eqe

With these six conditions, we can now state the  universality result for the Poisson case.

\begin{thm}[Poisson]\label{thm:univ_poiss}
Suppose that $\cH_\nn$ satisfies the assumptions in \eqref{eqn:assum_translate}-\eqref{eqn:assum_stable}, and let $f\in \cD(\cH)$. 
Let $\cP_\nn\sim \poisp{\nn}{f}$,
then
\[
\lim_{\nn\to\infty}\frac1\nn\mean{\cH_\nn(\cP_\nn)} = \cH^*,
\]
where $\cH^*$ is limit of the reference model, as defined in \eqref{eqn:assum_linear}.
\end{thm}

For the binomial case we will require two additional assumptions.

\noindent{\bf Polynomial bound.} We assume that almost surely,
\eqb\label{eqn:assum_poly_1}
\cH_n(\cX_n) \le C_3 n^q,
\eqe
for some $C_3,q>0$. 

\noindent{\bf Bounded difference.}
We assume there exists an event $A_n$, such that
\eqb\label{eqn:assum_A_n}
    \prob{A_n^c} \le e^{-n^a},
\eqe
for some $a>0$, and almost surely for all $\Delta < n$,
\eqb\label{eqn:assum_poly_2}
|\cH_n(\cX_{n}) - \cH_n(\cX_{n\pm \Delta})|\ind_{A_n} \le C_4 \Delta n^b,
\eqe
for some $C_4>0$ and $b<1/4$.
With these additional conditions, we can prove a stronger result for the binomial process.

\begin{thm}[Binomial]\label{thm:univ_binom}
Suppose that $\cH_n$ satisfies the assumptions in \eqref{eqn:assum_translate}-\eqref{eqn:assum_poly_2}, and let $f\in \cD(\cH)$.
Let $\cX_n\sim\binmp{n}{f}$, 
then
\[
\lim_{n\to\infty}\frac1n\mean{\cH_n(\cX_n)} = \cH^*,
\]
where $\cH^*$ is the limit of the reference model, defined in \eqref{eqn:assum_linear}. Furthermore, we have almost surely that
\[
\lim_{n\to\infty}\frac1n\cH_n(\cX_n) = \cH^*.
\]
\end{thm}

\subsection{Proofs}

We will focus our efforts on proving universality for the Poisson case. The proof for the binomial will follow from that.

\begin{proof}[Proof for Theorem \ref{thm:univ_poiss}]
The proof is divided into three main steps: (1) piecewise-constant densities, (2) bounded support densities, and (3) the entire class $\cD(\cH)$.

\noindent\underline{Step 1 -- piecewise constant densities:}\\
This step is the key to understand the roots of  universality.
We consider piecewise constant densities supported on $Q^d$. 
Fix $m>0$, and divide $Q^d$ into $M=m^d$ boxes of side-length $1/m$, denoted $Q_1,\ldots,Q_M$. Let $f$ be of the form
\[
f(x) = \sum_{i=1}^M c_iM\ind_{Q_i}(x),
\]
where $0\le c_i$ for all $i$,  and $\sum_i c_i = 1$. Further, suppose that the boxes are ordered so that the first $M^+$ boxes ($1\le M^+\le M$) have $c_i>0$, and the rest have $c_i=0$.

Let $\cP_\nn\sim\poisp{\nn}{f}$, and let $\cP_{\nn,i} = \cP_\nn\cap Q_i$. If $c_i=0$, then $\cP_{\nn,i} = \emptyset$ and $\cH_\nn(\cP_{\nn,i})=0$. 
Otherwise, note that $\cP_{\nn,i} \sim \poisp{ {c_i}\nn}{M\ind_{Q_i}}$.
Therefore, we have
\[
\lim_{\nn\to\infty}\frac1\nn\mean{\cH_\nn(\cP_\nn)} =\lim_{\nn\to\infty}\frac1\nn
\sum_{i=1}^M\mean{\cH_\nn(\cP_{\nn,i})} =
\frac1M\sum_{i=1}^{M_+}c_i\cH^* =  \cH^*,
\]
The first equality uses assumption \eqref{eqn:assum_additive}, and the second equality holds due to assumptions \eqref{eqn:assum_translate}-\eqref{eqn:assum_homog}.
This completes the proof for step 1.

In other words, the source of universality in this case is that in each of the small boxes $Q_i$, we have a uniform distribution, and by scale invariance, the limit will be $\cH^*$  as in the reference model (uniform in $Q^d$). Thus, taking a convex combination of the limits yields $\cH^*$ as well.

\noindent\underline{Step 2 -- bounded support densities:}\\
Here, we will prove universality for any density $f$ whose support is contained in $Q^d$. The idea is to approximate $f$ with a piecewise constant density, and bound the difference. Once proven for $Q^d$, by the shift and scale invariance, we can conclude that the same holds for any compactly supported distribution.

Let $f\in\cD(\cH)$ be a density with a bounded support, and without loss of generality, suppose that its support  is contained in $Q^d$.
For any $m\ge 1$, we denote $Q_{m,1},\ldots, Q_{m,M}$ the division of $Q^d$ into $M=m^d$ sub-cubes of side-length $1/m$. For every $1\le i \le M$, define
\[
c_{m,i} = \min_{x\in Q_{m,i}} f(x),
\]
and
\[
 f_m(x) = \sum_{i=1}^M c_{m,i}\ind_{Q_{m,i}}(x).
\]
Note that the sequence $f_1,f_2,f_4,f_8\ldots$, satisfies the conditions of the continuity assumption \eqref{eqn:assum_stable}. Therefore, for every $\eps>0$ we can find $m=2^\ell$ large enough, so that
\[
\lim_{\nn\to\infty} \frac1\nn\abs{\mean{\cH_\nn(\cP_\nn)}-\meanx{\cH_\nn(\cP_{\nn,m})}} \le \eps.
\]
Since $\cP_{\nn,m}\sim\poisp{(1-\delta_m)\nn}{\tilde f_m}$ has piecewise constant density, we conclude from the previous step and \eqref{eqn:assum_homog} that
\[
    \lim_{\nn\to\infty}\frac1\nn\meanx{\cH_\nn(\cP_{\nn,m})} = (1-\delta_m)\cH^*.
\]
Therefore,
\[
(1-\delta_m)\cH^* -\eps\le \liminf_{\nn\to\infty} \frac1\nn\mean{\cH_\nn(\cP_\nn)} \le\limsup_{\nn\to\infty} \frac1\nn\mean{\cH_\nn(\cP_\nn)} \le (1-\delta_m)\cH^*+\eps.
\]
As $\eps$ and $\delta_m$ can be arbitrarily small, this completes the proof for step 2.

\noindent\underline{Step 3 -- unbounded support:} \\
Let $f\in \cD(\cH)$ be an arbitrary density with a non-compact support.
Take a sequence of radii  $R_1< R_2 < R_3<\ldots < \infty$ and define
 \[ f_i(x) = f(x)\ind\set{|x|\le R_i},
 \]
 then the sequence $f_1,f_2,\ldots,$ satisfies the conditions of assumption \eqref{eqn:assum_stable}.
 Since $f_i$ has a compact support, we can conclude from the previous step and \eqref{eqn:assum_homog} that
 \[
    \lim_{\nn\to\infty}\frac1\nn \mean{\cH_\nn(\cP_{\nn,i})} = (1-\delta_i)\cH^*.
 \]
 Similarly to step 2, we can take $\eps$ and $\delta_i$ small enough, and use \eqref{eqn:assum_stable}  to complete the proof.
\end{proof}

Next, we prove universality for the binomial process, including a strong law of large numbers. The proof will rely on a  generic martingale-based statement for the law of large numbers that we present in the appendix,  Lemma \ref{lem:exp_bounds}.

\begin{proof}[Proof of Theorem \ref{thm:univ_binom}]
Consider the Poisson process $\cP_\nn$, with $\nn = n$, and let $N=|\cP_n|$. Then
\[
\mean{\cH_n(\cP_{n})}  
= \meanx{\cH_n(\cP_{n})\ ;\ {|N-{n}| > n^{3/4}/2}} + \meanx{\cH_n(\cP_{n})\ ;\ {|N-n| \le n^{3/4}/2}}.
\]
Using the Cauchy-Schwartz inequality, we have
\[
\Delta_1(n) :=\meanx{\cH_n(\cP_{n})\ ; \ {|N-n| > n^{3/4}/2}} \le \param{\meanx{\cH_n^2(\cP_{n})}\probx{|N-n| > n^{3/4}/2}}^{1/2}.
\]
Using the concentration of the Poisson distribution (e.g., in \cite{penrose_random_2003}) and \eqref{eqn:assum_poly_1}, we have
\eqb\label{eqn:mean_h_N}
\Delta_1(n) \le \sqrt{2} n^{q} e^{-n^{1/2}/18}=  o(n).
\eqe
Next, define 
\[
\splitb
\Delta_2(n) &:= \meanx{\cH_n(\cP_n)-\cH_n(\cX_n)\ ; \ |N-n|\le n^{3/4}/2} \\
&=\sum_{m=n_1}^{n_2}( \mean{\cH_n(\cX_m)} -\mean{\cH_n(\cX_n)})\prob{N=m},
\splite
\]
where $n_1=\ceil{n-n^{3/4}/2}$, $n_2 = \floor{n+n^{3/4}/2}$.
Recall the definition of the event $A_n$ in \eqref{eqn:assum_A_n}-\eqref{eqn:assum_poly_2}.
Then for all $n_1\le m \le n_2$  
\[
\splitb
\abs{\mean{\cH_n(\cX_m)}-\mean{\cH_n(\cX_n)}} &\le \mean{\abs{\cH_n(\cX_m)-\cH_n(\cX_n)}\ind_{A_n}}+ \mean{\abs{\cH_n(\cX_m)-\cH_n(\cX_n)}\ind_{A_n^c}}\\
&\le C_4|n-m|n^{b} + n_2^{q} e^{-n^a},\\
&\le C_4 n^{3/4+b} + n_2^{q} e^{-n^a},
\splite
\]
Therefore, if $b<1/4$, we have $\Delta_2(n) = o(n)$.
Combined with \eqref{eqn:mean_h_N}, we have
\[
\mean{\cH_n(\cX_n)}\probx{|N-n|\le n^{3/4}/2} + \Delta_2(n) \le \mean{\cH_n(\cP_n)} \le \mean{\cH_n(\cX_n)} + \Delta_1(n) + \Delta_2(n).
\]
Dividing by $n$ and taking the limit we have that 
\[
    \limninf \frac1n\mean{\cH_n(\cX_n)} = \cH^*.
\]
In addition, \eqref{eqn:assum_poly_1}-\eqref{eqn:assum_poly_2} imply that Lemma \ref{lem:exp_bounds} holds, and therefore for all $\eps>0$, we have
 \[
        \prob{\abs{\cH_n(\cX_n)-\mean{\cH_n(\cX_n)}} \ge \eps n} \le  e^{-n^\gamma}.
 \]
 Since the right hand side is summable, using the Borel-Cantelli Lemma, we conclude that almost surely,
 \[
    \limninf \frac1n\mean{\cH_n(\cX_n)-\mean{\cH_n(\cX_n)}} = 0,
\]
which concludes the proof.
\end{proof}

\section{Critical faces}\label{sec:crit_faces}

With the general framework in hand, our goal is to analyze the ensemble of $k$-cycles formed in the \emph{entire} thermodynamic limit ($nr^d = \lambda$ for any $\lambda\in(0,\infty)$). In order to bound the tail distribution in the thermodynamic limit, we will use  the notion of \emph{critical faces}.
Consider either the \v Cech filtration ($\set{\cC_r}_{r=0}^\infty$) or the Vietoris-Rips filtration ($\set{\cR_r}_{r=0}^\infty$). As we increase the parameter $r$, new simplices enter the filtration, potentially changing the homology of the complex. We refer to the simplices that facilitate these changes in homology as \emph{critical faces}. 

For the \v Cech complex, the Nerve Lemma \cite{borsuk_imbedding_1948} implies that $\cC_r(\pointset)$ is homotopy equivalent to $B_r(\pointset) = \cup_{x\in\pointset} B_r(x)$, and the latter is the sub-level set of the distance function to $\pointset$. Using Morse theory for the distance function \cite{bobrowski_distance_2014}, we can  define a unique $k$-dimensional critical face for every new $k$-cycle generated, and a unique $(k+1)$-dimensional critical face for every $k$-cycle terminated (becoming trivial).
Furthermore, precise geometric conditions for $k$-faces to be critical are provided. The main reason we can do so is that the random distance function is a min-type function \cite{gershkovich_morse_1997} which behaves like a Morse function. This implies that at every radius $r$ there can be at most a single change in the homology. 

To identify the critical faces, let $\cX^{(k)}$ denote the set of all $k$-simplices in $\cX$ (i.e., all subsets of size $k+1$), and let $\cY\in\cX^{(k)}$. Assuming general position, the points in $\cY$ can be placed on a unique $(k-1)$-sphere, whose center and radius we denote by $c(\cY)$ and $\rho(\cY)$. Denote by $B(\cY)$ the open $d$-dimensional ball centered at $c(\cY)$ with radius $\rho(\cY)$, and by $\sigma(\cY)$ the open $k$-simplex spanned by $\cY$ (i.e., the convex hull of $\cY$). The $k$-face $\cY$ is then critical if only if: (a) $c(\cY)\in \sigma(\cY)$, and (b) $B(\cY)\cap \cX = \emptyset$. Thus, defining $\tilde F_k(r)$ as the number of critical $k$-faces entering the filtration in the interval $[0,r]$, we have 
\eqb\label{eqn:crit_cech}
\tilde F_k(r) = \sum_{\cY\in \cX^{(k)}} h(\cY)\indf{B(\cY)\cap \cX = \emptyset} \indf{\rho(\cY)\le r},
\eqe
where $h(\cY) := \indf{c(\cY)\in \sigma(\cY)}$.

For the Vietoris-Rips complex, the situation is more complicated. When a new edge enters the complex at radius $r$, multiple changes can occur in homology simultaneously, across different homological degrees. Furthermore, apart from the edge, it is impossible to uniquely assign a single simplex to any of the changes. Fortunately, we do not need to actually identify the critical faces, but merely to count them (or more precisely, the changes in homology).  Suppose that $x_1,x_2\in \cX$, and $|x_1-x_2|=\rho$. Define $\lk((x_1,x_2);\cX)$ to be the link of the edge $(x_1,x_2)$, in the complex generated by $\cX$, at radius $\rho$ (see definition in Section \ref{sec:topo}). Lemma \ref{lem:crit_rips} states that the number of critical $k$-faces that are added at $\rho$ is  equal to ${\beta}_{k-2}(\lk((x_1,x_2);\cX))$ for $k\geq 2$, where $\beta_k$ is the \emph{reduced Betti number}.  Therefore, for the Vietoris-Rips complex ($k\ge 2$) we have, 
\eqb\label{eqn:crit_rips}
\tilde F_k(r) = \sum_{\set{x_1,x_2}\in \cX^{(1)}} {\beta}_{k-2}(\lk((x_1,x_2); \cX)) \indf{|x_1-x_2| \le r}.
\eqe
%
In the random setting, 
we will use $\tilde F_k(r)$ for the case when $\nn r^d = \lambda$ (for Poisson, or $nr^d=\lambda$ for binomial). Furthermore, we will sometime study  ranges of values $(\lambda_1,\lambda_2]$. Thus, from here onwards we will use 
\[
\splitb
F_k(\lambda) &:= \tilde F_k((\lambda/\nn)^{1/d}),\\
F_k(\lambda_1,\lambda_2) &:= F_k(\lambda_2)- F_k(\lambda_1).
\splite
\]
In both the \v Cech and Vietoris-Rips complex, we can write
\[
F_k(\lambda) = F_k^+(\lambda) + F_k^-(\lambda),
\]
where $F_k^+(\lambda)$ counts critical $k$-faces that generate new $k$-cycles (\emph{positive} faces), and $F_k^-(\lambda)$ counts critical $k$-faces that terminate existing $(k-1)$-cycles (\emph{negative} faces).

\vspace{5pt}
\noindent{\bf A word on notation:} Some of the limits we prove in the following are \emph{universal}, in which case we use a star ($*$), to denote that the limit is the same as the reference model (also denoted with $*$), independently of the underlying distribution. Other limits, which are distribution-specific, will be denoted by a diamond ($\diamond$).

\begin{prop}\label{prop:crit_lambda}
Let $f:\R^d\to\R$ be bounded  density function, and let $\cP_\nn\sim\poisp{\nn}{f}$. Let $F_{k,\nn}(\lambda)$ be the number of critical $k$-faces for either the \v Cech or the Vietoris-Rips filtration, generated by $\cP_\nn$.
Then,
\[
\lim_{\nn\to\infty}\frac1\nn\mean{F_{k,\nn}(\lambda)} = F^\diamond_k(\lambda;f),
\]
where $F^\diamond_k(\lambda;f)>0$ is given in \eqref{eqn:C_lambda_cech} for the \v Cech, and \eqref{eqn:C_lambda_rips} for the Vietoris-Rips. 
Furthermore, in both cases we have
\[
\lim_{\lambda\to\infty} F_k^\diamond(\lambda;f) = F_k^*,
\]
where $F^*_k>0$ is given in \eqref{eqn:C_infty_cech} for the \v Cech, and \eqref{eqn:C_infty_rips} for the Vietoris-Rips.
Most notably,  $F_k^*$ is independent of $f$ (but does depend on $d,k$ and the filtration type).
\end{prop}

The result of Proposition \ref{prop:crit_lambda} holds more generally for  Poisson processes with intensity $\nn f$, even when $f$ is not a probability measure. This will be useful later.
\begin{cor}\label{cor:crit_lambda_gen}
    Let $f:\R^d\to\R$ be a bounded and integrable function. Define $c = \int_{\R^d} f(x)dx$, and $\tilde f = f/c$. Take $\cP_\nn \sim\poisp{c\nn}{\tilde f}$, then 
    \[
\lim_{\nn\to\infty}\frac1\nn \mean{F_{k,\nn}(\lambda)} =  F_k^\diamond(\lambda;f),
    \]
    and 
    \[     \lim_{\lambda\to\infty} F_k^\diamond(\lambda;f) = c F_k^*,
    \]
where $F_k^\diamond(\lambda;f)$ and $F_k^*$ are the same as in Proposition \ref{prop:crit_lambda}.
\end{cor}

The following statements aim to cover the critical faces in the \emph{entire} thermodynamic limit (as opposed to a fixed $\lambda$). The proofs require the notion of `good' density functions $f$ -- see   Section \ref{sec:prob_set} for the definition of this class of functions. 

\begin{lem}\label{lem:good}
Let $f:\R^d\to\R$ be a good density function. There exists $\lmax = \lmax(\nn)\to\infty$, so that 
\[
\lim_{\nn\to\infty} \frac1\nn\mean{F_{k,\nn}(\lmax)} = F^*_k = \lim_{\lambda\to\infty} F^\diamond_k(\lambda;f).
\]
\end{lem}
In other words, for good densities, we can choose a radius beyond the thermodynamic limit that guarantees we count all critical faces within the thermodynamic limit, while the number of  critical faces at larger radii is  negligible. Note that (a) the limit of $\mean{F_{k,\nn}(\lmax)}$ in this case is \emph{universal}, and (b) the same limit holds for $F_{k,\nn}(\Lambda)$ for any $\Lambda\le \lmax$ that satisfies $\Lambda\to\infty$.

In the special case where the support of the density $f$ is compact and convex, a stronger statement can be made, and will be useful later.
\begin{lem}\label{lem:comp_conv}
Let $f:\R^d\to\R$ be a density with a compact and convex support $\cS$, such that $\inf_{\cS} f>0$. Then,
\[
\lim_{\nn\to\infty} \frac1\nn\mean{F_{k,\nn}(\infty)} = F^*_k = \lim_{\lambda\to\infty} F^\diamond_k(\lambda;f).
\]
\end{lem}

\section{Topological lemmas}\label{sec:topo}

In this section, we provide  two deterministic topological statements that will play a key role in the probabilistic proofs for Theorems \ref{thm:main_poisson} and \ref{thm:main_binom}. The first is about the critical faces in the Vietoris-Rips complex, and the second is a bound on the potential changes to persistence diagrams.

\subsection{Critical faces for the Vietoris-Rips filtration}

Unlike for the \v Cech filtration,  for the Vietoris-Rips filtration we are unable to appeal to Morse theory to define critical simplices. Note that the Vietoris-Rips complex is a special type of a \emph{clique} or \emph{flag} complex.
For these complexes, 
critical edges are well-defined. However, as 
higher dimensional simplices may  enter at the same time, the notion of  critical faces is difficult to define. Nevertheless, in this  section we quantify  the changes in the homology, 
analogously to the changes induced by critical points in Morse theory. Note that in this section, we do not make use of any metric properties of the Vietoris-Rips complex and thus our statements are phrased in terms of general clique complexes. We begin with some basic constructions. 

Let $K$ be an abstract simplicial complex. The (closed) star of a simplex $\sigma\in K$ is defined as 
\[
\clStr_K(\sigma) = \Cl(\{\tau \in K : \sigma \preceq \tau\}),\]
i.e., the smallest simplicial complex made by all simplices that contain $\sigma$.
The \emph{link} of a simplex $\sigma$ is defined as 
$$\lk_K(\sigma) =  \{\tau \in \clStr_K(\sigma) : \sigma\cap \tau = \emptyset\}.$$
Let $G$ be a graph, and 
 define $K = K(G)$ as the clique complex generated by $G$ (i.e., the maximal simplicial complex whose $1$-dimensional skeleton is $G$). 
 Take an edge $e\not\in G$, and set  $K':=K(G\cup \set{e})$. Since $K'=K\cup \clStr_{K'}(e)$, we can relate the homology groups of $K$, $K'$, and $K\cap \clStr_{K'}(e)$ via the reduced Mayer-Vietoris long exact sequence. Using the fact that $\clStr_{K'}(e)$ is contractible and  has trivial homology, we obtain
\begin{equation}\label{eqn:mvles}
\cdots\rightarrow{\Hg}_{k+1}(K') \xrightarrow{\delta_{k+1}} {\Hg}_k(K\cap \clStr_{K'}(e)) \xrightarrow{i_k} {\Hg}_k(K)\xrightarrow{q_k}  {\Hg}_k(K') \rightarrow \cdots 
\end{equation}
Denoting,
\begin{align*}
F^+_k (e) &=  \rk(\ker i_{k-1}),\\
F^-_{k} (e) &=\rk( \im i_{k-1}),
\end{align*}
the next lemma shows that for clique complexes, $F_k^+(e)$ represents the positive changes (new $k$-cycles generated) that occur when $e$ is added, and $F_k^-$ represents the negative changes (existing $(k-1)$-cycles terminated). The total number of critical $k$-faces is then $F_k(e) := F_k^+(e) + F_k^-(e)$.

\begin{lem}\label{lem:crit_rips}
Let $G$ be a graph and $e\not\in G$. Let $K=K(G)$ and $K' = K(G\cup\{e\})$ be the corresponding clique complexes. Then,
\begin{align}
\splitb
\beta_k(K') &= \beta_k(K) + F^+_k(e) - F^-_{k+1}(e),\quad k>0 \label{eqn:vietoris_critical_face}\\
\beta_0(K') &= \beta_0(K) - F^-_1(e),
\splite
\end{align}
and
\begin{equation}\label{eqn:rips_morse}
\splitb
F_k(e) &:= F_k^+(e) + F_k^-(e)  = \beta_{k-2}(\lk_{K'}(e)),\quad k>1 \\
F_1(e) &:= F_1^+(e) + F_1^-(e)  = \indf{\lk_{K'}(e) = \emptyset}.
\splite
\end{equation}
\end{lem}

\begin{rem}
We note that the proof below covers a 
 stronger statement than \eqref{eqn:vietoris_critical_face}.  Rather than the changes in the reduced Betti numbers, we can  characterize the actual changes in the homology groups. That is, considering the map $q_k: \Hg_k(K) \rightarrow \Hg_k(K')$ (induced by inclusion), then 
    our proof shows that $\rk(\coker q_k)  = F^+_k(e)$ and $ \rk(\ker q_k) = F^-_{k+1}(e)$.
\end{rem}

\begin{proof}
From  \eqref{eqn:mvles}, by exactness we have that
\begin{equation}\label{eqn:ranks_basic}
\rk({\Hg}_k(K'))  = \rk(\ker i_{k-1})  + \rk(\coker i_{k}).
\end{equation}
Furthermore,
\begin{align*}
\rk(\coker i_k) &= \rk({\Hg}_k(K)) - \rk(\im i_k)= \beta_k(K) -  F^-_{k+1}(e).
\end{align*}
For $k>0$, 
substituting this into \eqref{eqn:ranks_basic} yields \eqref{eqn:vietoris_critical_face}. For $k=0$, we note that $F_0^+(e)=0$ ($\ker i_{-1}= 0$). 

It remains to prove \eqref{eqn:rips_morse}. 
By the rank-nullity theorem, 
\[
\rk({\Hg}_{k-1}(K\cap \clStr_{K'}(e)))  = \rk(\ker i_{k-1}) +\rk(\im i_{k-1})
= F^+_{k}(e) +F^-_{k}(e).
\]
Thus, we need to show that
 ${\Hg}_{k-1}(K\cap \clStr_{K'}(e))\cong {\Hg}_{k-2}(\lk_{K'}(e))$. 
Let $e=(u,v)$ and take a cover of $X:=K\cap \clStr_{K'}(e)$ by $\mathcal{U} :=\clStr_{K}(u)\cap X$ and $\mathcal{V}:=\clStr_{K}(v)\cap X$. That $\mathcal{U}$ and $\mathcal{V}$ cover $X$ follows immediately, since any simplex in $\clStr_{K'}(e)\cap X$ must be in either $\clStr_{K}(u)$, $\clStr_{K}(v)$, or both.  
Also note that, by definition, we have
$\mathcal{U}\cap \mathcal{V} = \lk_{K'}(e)$.
To complete the proof, we apply the reduced Mayer-Vietoris sequence again, obtaining
$$ \cdots \rightarrow {\Hg}_{k}(\lk_{K'}(e)) \rightarrow {\Hg}_{k}(\mathcal{U})\oplus {\Hg}_{k}(\mathcal{V})  \rightarrow  {\Hg}_{k}(K\cap \clStr_{K'}(e)) \rightarrow \Hg_{k-1}(\lk_{K'}(e))\rightarrow\cdots$$
Since both $\clStr_{K}(u)\cap X$ and $\clStr_{K}(v)\cap X$ are contractible (as  stars of vertices), it follows by exactness that 
\[
\Hg_k(K\cap \clStr_{K'}(e)) \cong {\Hg}_{k-1}(\lk_{K'}(e)),
\]
completing the proof for first part of \eqref{eqn:rips_morse}. 

For the second part of \eqref{eqn:rips_morse} ($k=1$),
if $\lk_{K'}(e)\ne \emptyset$ then $K\cap \clStr_{K'}(e)$ is connected, implying that both $\im i_0$ and $\ker i_0$ are trivial, and  
 $F_1^+(e)=F_1^-(e) = 0$, so the equation holds trivially.  On the other hand, if $\lk_{K'}(e)= \emptyset$, then $u$ and $v$ are not connected in $K\cap \clStr_{K'}(e)$  and so  $F_1^+ +F_1^- = \rk(\ker i_0) + \rk(\im i_0)=1$.
  
\end{proof}

\subsection{Stability}

The proofs in Section \ref{sec:proofs} focus on subsets of the persistence diagram where $\pi(p)\ge \alpha$, for some fixed $\alpha$. We call the corresponding cycles ``$\alpha$-persistent," and denote their count by $\Pi_k(\alpha;\cF)$.  
A key argument in proving  Theorem \ref{thm:main_poisson} is the following result, showing that $\Pi_k(\alpha;\cF)$ is stable with respect to adding a single simplex to the filtration. Abusing notation, we write $\cF\cup{\sigma}$, to denote adding $\sigma$ at some point along the filtration (in a consistent way).

\begin{lem}\label{lem:addone_simplex}
Let $\cF$ be a filtration and $\sigma\not\in\cF$ such that $\cF\cup\sigma$ is a valid filtation. Let $\Pi_k(\alpha;\cF)$ be the number of $\alpha$-persistent $k$-cycles in $\dgm_k(\cF)$. For any $\alpha \in(1,\infty)$, if $\dim(\sigma) = k,k+1$, then
\[
\Pi_k(\alpha;\cF)-1 \leq \Pi_k(\alpha;\cF\cup \sigma) \leq \Pi_k(\alpha;\cF)+1.
\]
Otherwise, if $\dim(\sigma) \ne k,k+1$ then
\[
 \Pi_k(\alpha;\cF\cup \sigma) = \Pi_k(\alpha;\cF).
\]
\end{lem}

\begin{rem}\label{rem:addone}
The proof of Lemma \ref{lem:addone_simplex} below shows that if $\dim(\sigma)=k$ then
the difference in $\Pi_k(\alpha;\cF)$ comes from a single $k$-cycle, generated by $\sigma$, and hence has an infinite death time. Thus, when we later consider  cycles with a bounded death time ($\rmax$), the only difference in $\Pi_k(\alpha;\cF)$ stems from the case $\dim(\sigma)=k+1$.
\end{rem}

We begin with some notation.
Let $\cF':=\cF\cup \sigma$, and let $\dgm_k(\cF)$ and $\dgm_k(\cF')$ denote the respective persistence diagrams  and $i: \cF \hookrightarrow \cF' $ denote the 
natural inclusion map.  Let  $w(\sigma)$ be the filtration value function (i.e., the time at which $\sigma$ joins the complex). 

Recall that the persistence diagram is a representation of the algebraic structure of the corresponding persistence module, namely each point in the diagram corresponds to a rank-one summand in the decomposition \eqref{eqn:ph_decomp}.  We use $\gamma$ to denote the persistent homology class corresponding to a single point in $\dgm(\cF)$.  
In this section, we assume that the filtration has been extended to a total order (i.e., that the filtrations values $w$ are unique). In cases where multiple simplices enter the filtration at the same time, the ties are broken in an arbitrary but consistent way (i.e., respecting the conditions for simplicial complexes).

Note that inserting a $k$-simplex into a simplicial complex either creates a new class in the $k$-th homology, in which case the simplex is referred to as a \emph{positive},  or it bounds an existing class in the $(k-1)$-th homology, and is referred to as \emph{negative} (c.f., \cite{delfinado_incremental_1993}). 
We define the function $\pairing_{\cF}(\gamma)$  as a matching between a class $\gamma$ and the negative simplex in $\cF$ that kills it.  This matching is well-defined due to the  assumption on the total order. 
\begin{defn}
Let $\cF$ be a filtration,  $\sigma\not \in \cF$, and $\cF' = \cF\cup\sigma$.
 The \emph{cascade} of $\sigma$, denoted    $\mathrm{Csc}_{\cF}(\sigma)$, is the set of points in the persistence diagram, whose values change when $\sigma$ is added to $\cF$, and is defined as,
$$\Csc_{\cF}(\sigma) = \set{(\birthx_j,\deathy_j) \in \dgm_k(\cF) : \pairing_{\cF}(\gamma_j) \neq \pairing_{\cF'}(i_*(\gamma_j))  }  $$
where $i_*$ is map on homology induced by the inclusion $i:\cF\to\cF'$. 
\end{defn}

The key technical lemma in this section is regarding the structure of a cascade.

\begin{lem}\label{lem:cascade}
Suppose that 
$\Csc_{\cF}(\sigma)= \{(\birthx_1,\deathy_1),\ldots, (\birthx_m,\deathy_m)\}$ where $\birthx_1\geq \birthx_2 \geq \ldots \geq \birthx_m$, and denote the corresponding classes by $\gamma_1,\ldots\gamma_m$. If $\sigma$ is a negative simplex in $\cF'$, then
\begin{enumerate}
\item[(A)] $w(\sigma)\leq \deathy_1\leq \deathy_2 \leq\cdots \leq \deathy_m$,
\item[(B)] $\pairing_{\cF'}(i_*(\gamma_1)) = \sigma$,
\item[(C)] $\pairing_{\cF'}(i_*(\gamma_j)) = \pairing_{\cF}(\gamma_{j-1})$ for $j>1$.
\end{enumerate}
\end{lem}

Note that from $(A)$ we can conclude that $(\birthx_1,\deathy_1) \subset  (\birthx_2,\deathy_2) \subset \cdots \subset(\birthx_m,\deathy_m)$, justifying the name `cascade'.
We postpone proving this lemma to later, but we use it next to prove the main result. 

\begin{proof}[Proof of Lemma \ref{lem:addone_simplex}]
First, note that if $\dim(\sigma)\ne k,k+1$, then $\sigma$ has no effect on $\Hg_k$. Next, suppose that $\dim(\sigma) = k$. In this case there are two possible cases which can affect $\Pi_k$:
\begin{itemize}
    \item $\sigma$ is   a positive simplex in $\cF'$,
    \item $\sigma$ is a negative simplex in $\cF'$, and there is another $k$-simplex $\sigma'\in \cF$ that is negative in $\cF$ and positive in $\cF'$ (i.e., $w(\sigma) < w(\sigma')$ with 
    $\sigma$ killing the $(k-1)$-cycle that $\sigma'$ used to)\footnote{More precisely, this occurs when the chains $\partial\sigma'$ and $\partial\sigma$ are  homologous, 
    and when $\sigma$  enters both become a boundary. Thus, the insertion of  $\sigma'$ forms a new cycle in $\cF'$.}.
\end{itemize}
In both cases, a new $k$-th homology class is generated whose cycle representative includes $\sigma$. This new class must be  infinite in $\cF'$, since by assumption $\sigma$ has no cofaces in $\cF'$. Hence, we can deduce that
$\dgm_k(\cF') = \dgm_k(\cF) \cup \{(w(\sigma), \infty)\}. $
That is, we add one point to the diagram with birth at $w(\sigma)$ and infinite death time. Hence, $\Pi_k(\alpha;\cF') = \Pi_k(\alpha;\cF)+1$ for all $\alpha$. 

The last case is when 
 $\sigma$ is a negative $ (k+1)$-simplex in $\cF'$. Here, by definition, all the points in the diagram which contribute to $\Delta \Pi_k(\alpha) := |\Pi_k(\alpha;\cF)-\Pi_k(\alpha;\cF')|$ are in 
$\Csc_{\cF}(\sigma)$.  Indexing the points as in Lemma \ref{lem:cascade}, we first observe that since the $\birthx_j$'s are decreasing and the  $\deathy_j$'s are increasing ((A)~in Lemma \ref{lem:cascade}), the ratio ${\deathy_j}/{\birthx_j}$ is increasing in $j$. 
Denoting $w(\pairing_{\cF'}(i_*(\gamma_j)))$ by $\deathy'_j$, we observe that for $j>1$, $\deathy'_j = \deathy_{j-1} \leq \deathy_j$ and for $j=1$, $\deathy'_1 = w(\sigma) \leq \deathy_1$ ((B)~and (C)~in Lemma~\ref{lem:cascade}).
This implies that for all $j$, we have ${\deathy_j}/{\birthx_j}\geq {\deathy'_j}/{\birthx_j}$.

Let $j_0$ be the first point in $\Csc_{\cF}(\sigma)$ such that 
${\deathy_{j_0}}/{\birthx_{j_0}}\geq \alpha$.
We argue that for all $j\ne j_0$, the point $(\birthx_j,\deathy_j)\in \Csc_{\cF}(\sigma)$ does not contribute to $\Delta\Pi_k(\alpha)$.
If $j<j_0$, then $\deathy'_j/\birthx_j \le \deathy_j/\birthx_j < \alpha$, and therefore the $j$-th point does not contribute to  $\Delta\Pi_k(\alpha)$. If $j>j_0$, then $\deathy'_j = \deathy_{j-1}\ge \deathy_{j_0}$, and $\birthx_{j}\le \birthx_{j_0}$. Therefore, $\deathy_j/\birthx_j\ge \deathy'_j/\birthx_j \ge \deathy_{j_0}/\birthx_{j_0}\ge \alpha$, so here as well the $j$-th point does not contribute to $\Delta\Pi_k(\alpha)$. 

To summarize, out of all the points in $\dgm_k(\cF)$ that change as a result of adding $\sigma$, only a single point $(\birthx_{j_0},\deathy_{j_0})\in \Csc_{\cF}(\sigma)$ might contribute to $\Delta\Pi_k(\alpha)$, concluding the proof. 
\end{proof}

We conclude this section  by proving Lemma~\ref{lem:cascade}. The proof  relies on the following result from \cite{morozov_homological_2008}, describing  the possible changes in persistence diagrams in response to changing the order of the simplices in the filtration. As before, we assume a total ordering on the simplices in the filtration.
Somewhat surprisingly, swapping the order of two consecutive simplexes, can affect at most two points in the persistence diagram (see the Switch Lemma in 
\cite{morozov_homological_2008}).
The types of possible changes are listed in the next theorem. To describe these changes, a useful view for points in a persistence diagram is as pairings
between positive and negative simplices (along with their filtration values), which we refer to as the \emph{birth-death pairings}. Suppose that $(\birthx_1,\deathy_1),(\birthx_2,\deathy_2)\in \dgm_k(\cF)$. We say that the corresponding pairings are \emph{nested} if $(\birthx_1,\deathy_1) \subset (\birthx_2,\deathy_2)$ (as intervals in $\R$), and similarly we say the pairings are \emph{disjoint} if $(\birthx_1,\deathy_1) \cap(\birthx_2,\deathy_2) = \emptyset$.

\begin{thm}[Pairing Change Theorem, \cite{morozov_homological_2008}]\label{thm:pairing_change}
After switching the order of two consecutive simplices, the 
birth-death pairings may change only if the dimension of the
swapped simplices is the same. The possible  changes are:
\begin{enumerate}
\item[(1)] Two nested pairings swap their birth-simplices to remain nested.
\item[(2)] Two nested pairings swap their death-simplices to remain nested. 
\item[(3)] Two disjoint pairings swap a birth-simplex with a death-simplex to remain disjoint.  
\end{enumerate}
\end{thm}

\begin{figure}
\centering
\includegraphics[width=0.8\textwidth,page=8]{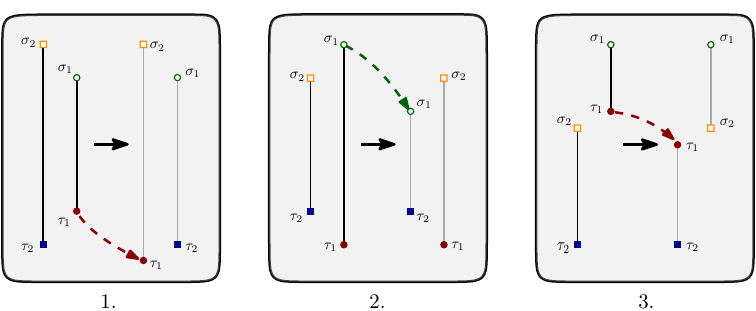}
\caption{\label{fig:switches} The three possible cases of pairing switches. We depict the two affected pairings, before and after the switch, where the height of the bottom points represent birth, and the height of the top ones represent death. The arrow indicates which simplex is moving backwards in the filtration and how the pairing changes. Note that moving simplices forward is simply the reverse of one of the four cases above.}
\end{figure}

We illustrate the possible changes in Figure~\ref{fig:switches}. In case (1),  two positive simplices are swapped (affecting the birth times). In case (2), two negative simplices are swapped (affecting death times). Finally, in case (3), a positive and negative simplex are swapped. We note that although the changes in the birth-death pairings is effectively the same in cases (1) and (2) (i.e., we swap the matching of the birth and death simplices),  
the difference is in which of the values changes: in case (1), the birth time changes value, while in cases (2), the death time changes value. 

\begin{proof}[Proof of Lemma \ref{lem:cascade}]
We prove the characterization inductively.
We first construct a filtration $\cF_0$, where $\sigma$ is appended to the end of the filtration $\cF$. We then construct a sequence of filtrations, $\cF_\ell$ such that at each step, $\sigma$ is swapped with an earlier simplex towards its final position in $\cF'$. We use Theorem~\ref{thm:pairing_change}, to show that the cascade satisfies the  properties (A)-(C) stated in Lemma \ref{lem:cascade} at each $\cF_\ell$ and refer to the intermediate filtration functions as $w_\ell$.

If $\sigma$ appended to the end of $\cF$ is positive, the cascade is empty and so the properties holds vacuously. By assumption, $\sigma$ at its final location in the filtration is negative, so there must be a step $\ell_0\ge 0$ where $\sigma$ becomes negative. 
This is case (3) in Theorem~\ref{thm:pairing_change}, where a positive and negative simplices switch. Hence, the cascade consists of a single point, whose death time corresponds to the simplex which was swapped with $\sigma$.  As there is only one point in the cascade,  (A) and (C) hold vacuously. As the death time of the point is given by the current filtration value of $\sigma$,  (B) holds as well. 

For steps $\ell>\ell_0$,
 we continue to move $\sigma$ backwards in the filtration. At each step, if no pairing switch occurs, then the cascade remains unchanged and so the (A)-(C) still hold. As cases (1a) and (1b) in Theorem~\ref{thm:pairing_change} are for swapping positive simplices, these changes cannot occur here. Likewise, case (3) cannot occur as $\sigma$ would need to be positive before the swap.

It remains to check the effect of case (2).
Suppose that $(\birthx_1,\deathy_1), (\birthx_2,\deathy_2) \in \dgm_k(\cF_{\ell-1})$ are the points affected by the swap, so that $w_{\ell-1}(\sigma) = d_2$ and $w_\ell(\sigma) = d_1'\le d_1$.
Since the pairings are nested before and after the swap, we have $\birthx_2\leq \birthx_1 \leq \deathy_1 \leq \deathy_2$, and 
$\birthx_2\leq \birthx_1 \leq \deathy'_1 \leq \deathy_1$. By the induction assumption, (A) holds before the swap, so ordering by decreasing values of $\birthx_j$, we have $d_2\leq d_j$ for $j>2$. 
After the swap,  $b_1\geq b_2$, 
and $w_\ell(\sigma)=\deathy'_1 \leq d_1$ and $d_1\leq d_3 \leq \cdots$ (recall that $d_2$ was changed into  $d'_1$). Reindexing by decreasing $b_j$, we obtain (A) again.

Let $\gamma_1$ and $\gamma_2$ denote the respective classes in $\PH_k(\cF_{\ell-1})$. Note that we can identify $\gamma_1$ with $\birthx_1$ and $\gamma_2$ with $\birthx_2$. Hence, after the swap, as the new point in the diagram is $(\birthx_1,\deathy'_1)$, we conclude that $\sigma$ is matched with  $\gamma_1$, which is precisely (B).
 As these are the only changes, at each step the properties hold, completing the proof. 
\end{proof}

\section{Proofs}\label{sec:proofs}

In this section we provide  the details needed for the proofs of Theorems \ref{thm:main_poisson} and \ref{thm:main_binom}. We first introduce some useful notation. We use the following asymptotic notation:
\begin{itemize}
    \item 
    $a_n \sim b_n$ if $a_n = \Theta(b_n)$.
    \item 
    $a_n\approx b_n$ if $a_n/b_n \to 1$.
    \item 
    $a_n \ll b_n$ if $a_n = o(b_n)$.
\end{itemize}
Additionally, for $f:\R^d\to \R$, and a Borel set  $A\subset\R^d$, we denote
$
    f(A) := \int_A f(x)dx.
$
We also denote $B_r(A) := \cup_{x\in A}B_r(x)$.
Finally, we use `$C$' to represent constants where the actual value plays no role in the result, and needs no tracking. Thus, the actual value of $C$ may change, even within the same equation.

\noindent{\bf The largest death time.} 
Recall that $\rmax = \rmax(n)$ (and correspondingly $\lmax = n\rmax^d$) is bounding the death times that we allow. While none of the results explicitly depend on the value of $\rmax$, the proof techniques impose different conditions in different settings. For Theorem \ref{thm:main_poisson}, if $f$ is Type I or II, we can take $\lmax =n^{\frac1{dk+2}}$, while for Type III it is sufficient that $\rmax = O\param{1/{\log n}}$ (or $\lmax = O\param{{n}/{(\log n)^d}}$). For the law of large numbers in Theorem \ref{thm:main_binom}, we can take $\lmax = n^{\frac{1}{4(k+2)}}$. In all these cases, the important part is that $\rmax\to 0$ while  $\lmax\to\infty$, so that we account for all cycles within the thermodynamic limit.

\subsection{Persistent cycles -- Poisson}

Our goal in this section is to prove Theorem \ref{thm:main_poisson} using Theorem \ref{thm:univ_poiss}. In this case, we fix $\alpha$, and take $\cH_\nn(\cP_\nn) = \Pi_{k,\nn}(\alpha)$. Note that the parameter $\nn$ controls both the rate of the Poisson process ($\cP_\nn$) and the threshold death value ($\rmax = \rmax(\nn)$). We start by analyzing the reference model, and then move on to the general case. We note that while the analysis of critical faces is different between the \v Cech and Vietoris-Rips complexes (see Sections \ref{sec:crit_faces} and  \ref{sec:crit_faces_proofs}), the proofs in this section are the same for both filtrations.

\subsubsection*{The reference model}

Recall that our reference model is $\cP_\nn = \cP_\nn^* \sim \poisp{\nn}{\ind_{Q^d}}$.
We start with the following lemma, proving the limit for $\bar\Pi_{k,\nn}(\R)$ (expected number of \emph{all} points in the persistence diagram with $\death\le \rmax$), and thus the finiteness of the limiting measure.

\begin{lem}\label{lem:total_pers}
Let $f:\R^d\to\R$ be a density with a compact and convex support $\cS$, such that $\inf_{\cS} f>0$. Let $\cP_\nn\sim\poisp{\nn}{f}$. Then,
\[
\lim_{\nn\to\infty} \frac1\nn \bar\Pi_{k,\nn}(\R)= \Pi_k^*(\R):= \sum_{j=0}^k (-1)^{k-j}F^*_j,
\]
where $F_j^*$ are the limits defined in Proposition \ref{prop:crit_lambda}, and $F^*_0=1$. Note that since $F^*_j$ are independent of $f$, so is $\Pi_k^*(\R)$.
\end{lem}

\begin{proof}
Consider either the \v Cech or the Vietoris-Rips filration. 
Let $\cX\subset \R^d$ be a finite set, and 
let $F_k$ be the total number of critical $k$-faces in the corresponding filtration (i.e., without any restriction on the radius). 
We can decompose $F_k$  into $F_k^++F_k^-$, standing for the 
{positive} and {negative} faces. For every $k\ge 1$ we have $\Pi_{k}(\R) = F_k^+ = F_{k+1}^-$ (since every cycle  born is eventually  killed\footnote{Importantly,  this holds in this case only because we do not threshold the death-times with $\rmax$.}).
Finally, the total number of components generated is $|\cX|$, and all but one of them die. Therefore, 
\[
\splitb
F_1^- &=  F_0 := |\cX|-1,\\
F_2^- &= F_1^+ = F_1 - F_0,\\
F_3^- &= F_2^+ = F_2 - F_1 + F_0,\\
\splite
\]
and by induction, we have
\eqb\label{eqn:Pik_morse}
\Pi_k(\R) = \sum_{j=0}^k (-1)^{k-j} F_j.
\eqe
Switching  to the (random) reference model, recall that $\Pi_{k,\nn}$ counts cycles with death bounded by $\rmax$. Define $\widetilde\Pi_{k,\nn}$ the total count without the $\rmax$ restriction. Then by bounding the negative faces, 
\[
\widetilde\Pi_{k,\nn}(\R)-\Pi_{k,\nn}(\R)\le F_{k+1,\nn}(\lmax,\infty),
\]
and  using Lemma \ref{lem:comp_conv} we have,
\[
\lim_{\nn\to\infty}\frac1\nn \Pi_{k,\nn}(\R) = \lim_{\nn\to\infty}\frac1\nn\widetilde\Pi_{k,\nn}(\R).
\]
Finally, applying \eqref{eqn:Pik_morse} to $\widetilde\Pi_{k,\nn}(\R)$, 
noting that for $\cP_\nn$ we have $\mean{F_{0,\nn}(\infty)} = \nn-1$, and using Lemma \ref{lem:comp_conv}, completes the proof.
\end{proof}

Next, we show that the limits for $\alpha$-persistent $k$-cycles in   the reference model exist.
\begin{lem}\label{lem:pers_assum_linear}
Fix $\alpha \in (1,\infty)$ and let $\Pi_{k,\nn}^*(\alpha) := \Pi_{k,\nn}^*([\alpha,\infty))$ be the number of $\alpha$-persistent $k$-cycles generated by $\cP_\nn^*\sim \poisp{\nn}{\ind_{Q^d}}$. Then
\[
\lim_{\nn\to\infty}\frac1\nn \meanx{\Pi_{k,\nn}^*(\alpha)} = \Pi^*_k(\alpha),
\]
for some $\Pi_k^*(\alpha)>0$.
\end{lem}

\begin{rem}
Note that the case $\alpha=1$ was covered in Lemma \ref{lem:total_pers} (where the limit is $\Pi_k^*(\R)$). Ideally (and supported by experiments), we should have $\lim_{\alpha\to 1}\Pi_k^*(\alpha) = \Pi_k^*(\R)$. However, the proof for this claim remains as future work.
\end{rem}

\begin{proof}
Fix any $\alpha > 1$. We will decompose $\Pi_{k,\nn}^*(\alpha)$ into different radii ranges. Let $I\subset [0,\infty)$, and define $\Pi_{k,\nn}^*(\alpha; I)$ as the number of points in $\dgm_k(\cP_\nn^*)$ with $\pi(p)\ge \alpha$, and with $\death(p)=\rho$ such that $\nn \rho^d \in I$. For any $\lambda>0$, we can write
\eqb\label{eqn:decomp_lambda}
 \Pi_{k,\nn}^*(\alpha) =  \Pi_{k,\nn}^*(\alpha;(0,1/\lambda])+ \Pi_{k,\nn}^*(\alpha; (1/\lambda,\lambda]) + \Pi_{k,\nn}^*(\alpha; (\lambda,\lmax)).
\eqe
To bound the first and last term, it is sufficient to bound the number of negative critical $(k+1)$-faces appearing in the corresponding range. We will do so using Proposition \ref{prop:crit_lambda}. Whenever it is clear what $f$ is (e.g., $f=\ind_Q$ in this case), we will use $F_k^\diamond(\lambda) :=F_k^\diamond(\lambda;f)$.
For the first term, we have
\[
\meanx{\Pi_{k,\nn}^*(\alpha;(0,1/\lambda])} \le \meanx{F_{k+1,\nn}(1/\lambda)} \approx \nn F^\diamond_{k+1}(1/\lambda).
\]
Since $\lim_{\lambda\to \infty} F_{k+1}^\diamond(1/\lambda) = 0$, for large enough $\lambda$ we have  $F^\diamond_{k+1}(1/\lambda) \le \eps/3$. 
For the last term in \eqref{eqn:decomp_lambda}, we use Lemma \ref{lem:good}
\[
\meanx{\Pi_{k,\nn}^*(\alpha;(\lambda,\lmax))} \le \meanx{F_{k+1,\nn}(\lambda,\lmax)} \approx \nn(F^*_{k+1}-F^\diamond_{k+1}(\lambda)).
\]
By the definition of $F_{k+1}^*$ (Proposition \ref{prop:crit_lambda}), for large enough $\lambda$ we have $(F^*_{k+1}-F^\diamond_{k+1}(\lambda))\le\eps/3$.

Finally, from \cite[Theorem 1.5]{hiraoka_limit_2018}, we know that for $\cP_\nn^*$, and for any fixed $\lambda$, we have
\[
    \lim_{\nn\to\infty} \frac1\nn \meanx{\Pi_{k,\nn}^*(\alpha;(1/\lambda,\lambda])} = \cL(\alpha;\lambda),
\]
for some $\cL(\alpha;\lambda)>0$. Note that $\cL(\alpha;\lambda)$ is increasing in $\lambda$, and is also bounded by $\Pi_k^*(\R)$, and therefore there exists a limit $\cL(\alpha) = \lim_{\lambda\to\infty}\cL(\alpha;\lambda)$. Thus, for large enough $\lambda$, we have $|\cL(\alpha;\lambda)-\cL(\alpha)|<\eps/3$.
Combining the three bounds, we have
\[
\cL(\alpha)-\eps/3\le \liminf_{\nn\to\infty} \frac1\nn\meanx{ \Pi_{k,\nn}^*(\alpha)} \le  
 \limsup_{\nn\to\infty} \frac1\nn\meanx{ \Pi_{k,\nn}^*(\alpha)}
\le \cL(\alpha)+ \eps.
\]
As this holds for any $\eps>0$, setting $\Pi_k^*(\alpha) = \cL(\alpha)$ concludes the proof.

\end{proof}

\subsubsection*{General distributions}

In this section we  prove universality for the class of good densities (see Section \ref{sec:prob_set}).

\begin{lem}\label{lem:limit_good}
    Let $f$ be a good density, and $\cP_\nn\sim \poisp{\nn}{f}$. Then for all $\alpha \ge 1$, we have
    \[
    \lim_{\nn\to\infty}\frac1\nn \mean{\Pi_{k,\nn}(\alpha)} = \Pi_k^*(\alpha)
    \]
    where the limits $\Pi_k^*(\alpha)$ are the same as for the reference model, presented in Lemmas \ref{lem:total_pers} and \ref{lem:pers_assum_linear}, respectively.
\end{lem}

\begin{proof}
     Fix $\alpha\ge 1$. If we can show that the conditions of Theorem \ref{thm:univ_poiss}  apply here (i.e., with $\cH_\nn(\cP_\nn)=\Pi_{k,\nn}(\alpha)$) that will prove the lemma. The translation-invariance condition \eqref{eqn:assum_translate} holds immediately for $\Pi_{k,\nn}(\alpha)$. Assumption \eqref{eqn:assum_linear} is proved by Lemma \ref{lem:pers_assum_linear}. Finally, Lemmas \ref{lem:pers_assum_scale}-\ref{lem:pers_assum_stable} below, verify that conditions \eqref{eqn:assum_scale}-\eqref{eqn:assum_stable} hold as well. 

\end{proof}

We are therefore left with proving the following lemmas. 

\begin{lem}\label{lem:pers_assum_scale}
Assumption \eqref{eqn:assum_scale} holds for $\Pi_{k,\nn}(\alpha)$.
\end{lem}

\begin{lem}\label{lem:pers_assum_homog}
Assumption \eqref{eqn:assum_homog} holds for $\Pi_{k,\nn}(\alpha)$.
\end{lem}

\begin{lem}\label{lem:pers_assum_additive}
Assumption \eqref{eqn:assum_additive} holds for $\Pi_{k,\nn}(\alpha)$.
\end{lem}

\begin{lem}\label{lem:pers_assum_stable}
Assumption \eqref{eqn:assum_stable} holds for $\Pi_{k,\nn}(\alpha)$.
\end{lem}

\begin{proof}[Proof of Lemma \ref{lem:pers_assum_scale}]
Fix $\alpha \ge 1$, and let $\Pi_{k,\nn}^*(\alpha)$ and $\Pi_{k,\nn}^{(\eps)}(\alpha)$ be the persistence measures corresponding to the reference model $\cP_\nn^*\sim\poisp{\nn}{\ind_{Q^d}}$ and its scaled version $\cP_\nn^{(\eps)}\sim\poisp{\nn}{\eps^{-d}\ind_{\eps Q^d}}$. We need to  show that
\[
\lim_{\nn\to\infty}\frac1\nn{\meanx{\Pi_{k,\nn}^{(\eps)}(\alpha)}} = \lim_{\nn\to\infty}\frac1\nn{\meanx{\Pi_{k,\nn}^*(\alpha)}}.
\]
Recall that $\Pi_{k,\nn}(\alpha)$ counts $k$-cycles with $\pi(p) \ge\alpha$, and $\death(p)\le \rmax$. If we define $\widetilde\Pi_{k,\nn}$ similarly, but without the restriction on the death-time, then $\widetilde\Pi^{(\eps)}_{k,\nn}(\alpha)$ and $\widetilde\Pi^{*}_{k,\nn}(\alpha)$ have the same distribution, and in particular the same expected value. Therefore, it suffices to show that the difference ($\widetilde\Pi_{k,\nn}(\alpha)-\Pi_{k,\nn}(\alpha)$) is negligible.
To this end, note that \[
\widetilde\Pi_{k,\nn}(\alpha)-\Pi_{k,\nn}(\alpha) \le {F}_{k+1,\nn}(\lmax,\infty).
\]
Using Lemmas \ref{lem:good} and \ref{lem:comp_conv}, we have that
\[
\lim_{\nn\to\infty}\frac1\nn \param{\widetilde\Pi^*_{k,\nn}(\alpha)-\Pi^*_{k,\nn}(\alpha)}  = \lim_{\nn\to\infty}\frac1\nn \param{\widetilde\Pi^{(\eps)}_{k,\nn}(\alpha)-\Pi^{(\eps)}_{k,\nn}(\alpha)} = 0,
\]
completing the proof.
\end{proof}

\begin{proof}[Proof of Lemma \ref{lem:pers_assum_homog}]
Denote $\Pi_{k,\nn, c\nn}(\alpha)$ the number of $\alpha$-persistent $k$-cycles for $\cP_\nn$, with death time bounded by $\rmax(c\nn) < \rmax(\nn)$. We need to show that $\lim_{\nn\to\infty}{\mean{\Pi_{k,\nn,c\nn}}}/{\mean{\Pi_{k,\nn}}} = 1$.
Note that
\[ 
\Pi_{k,\nn}(\alpha) = \Pi_{k,\nn,c\nn}(\alpha)+ \Pi_{k,\nn}(\alpha; (\lmax(c),\lmax]),
\]
where $\lmax(c) = n\rmax^d(c\nn)$. Since $\rmax(c\nn) < \rmax(\nn)$, and $\rmax(c\nn)\to\infty$, using Lemma \ref{lem:good},
\[
\lim_{\nn\to\infty} \frac{\Pi_{k,\nn}(\alpha; (\lmax(c),\lmax])}{\Pi_{k,\nn(\alpha)}} = 0,
\]
 concluding the proof.
\end{proof}

\begin{proof}[Proof of Lemma \ref{lem:pers_assum_additive}]
Take any $M=m^d$, and define $\cP_\nn = \bigcup_i \cP_{\nn,i}$, where $\cP_{\nn,i}\sim\poisp{c_i\nn}{M\ind_{Q_i}}$ are independent. For all $1\le i \le M$, let $\cK_i$ denote the filtration (either \v Cech or Vietoris-Rips) generated by $\cP_{\nn,i}$, and denote by $\cK$ the filtration generated by $\cP_\nn$.
Defining $\Pi_{k,\nn}^{(i)}(\alpha)$ as the persistence measure for $\dgm_k(\cK_i)$, our goal is to prove
\eqb\label{eqn:pi_additive}
\lim_{\nn\to\infty} \frac1\nn \abs{\mean{\Pi_{k,\nn}(\alpha)} - \sum_{i=1}^M \meanx{\Pi_{k,\nn}^{(i)}(\alpha)}} = 0.
\eqe
Similarly to the proof of Lemma \ref{lem:pers_assum_linear}, for any $I\subset[0,\infty)$, define $\Pi_{k,\nn}^{(i)}(\alpha;I)$, as the number points $p\in\dgm_k(\cK_i)$ with $\pi(p)\ge\alpha$ and with $\death(p)=\rho$, such that  $\nn\rho^d\in I$. We will split the proof of \eqref{eqn:pi_additive} into $I=(0,\lambda]$, and $I=(\lambda,\lmax]$.

We start with $I=(0,\lambda]$. Fix $\lambda>0$, and 
let $\cK_i^{(\lambda)}\subset \cK_i$ be the filtration generated by $\cP_{\nn,i}$,  with simplices that appear at radius $\rho$ such that $\nn\rho^d \le \lambda$. Similarly, define $\cK^{(\lambda)}\subset \cK$ as the filtration generated by the full process $\cP_\nn$. Then
\[
\cK^{(\lambda)} = (\cK_1^{(\lambda)} \sqcup \cdots\sqcup \cK_M^{(\lambda)})\sqcup \Sigma,
\]
where $\Sigma$ contains all simplices generated by points from two or more different boxes, and with radius $\rho$ satisfying $\nn\rho^d\le \lambda$. 
Note that $\cK_1^{(\lambda)},\ldots,\cK_M^{(\lambda)}$ are disjoint filtrations, and therefore
\[
\dgm_k(\cK_1^{(\lambda)}\sqcup\cdots\sqcup\cK_M^{(\lambda)}) = \dgm_k(\cK_1^{(\lambda)})\sqcup\cdots\sqcup \dgm_k(\cK_M^{(\lambda)}).
\]
Thus, using Lemma \ref{lem:addone_simplex} and Remark \ref{rem:addone}, we have
\[
\abs{\Pi_{k,\nn}(\alpha;(0,\lambda]) - \sum_i \Pi_{k,\nn}^{(i)}(\alpha;(0,\lambda])} \le \Sigma_{k+1},
\]
where $\Sigma_{k+1}$ counts  all $(k+1)$-simplices in $\Sigma$. Note that  these simplices must contain points from several boxes, and their radius is bounded by $r = (\lambda/\nn)^{1/d}$. Thus,  they are necessarily generated by points lying in $\partial_r := B_r\param{\bigcup_i \partial Q_i}$, where $\partial Q_i$ is the boundary of $Q_i$. In addition, since $Ar\le \vol(\partial_r) \le Br$ for some $A,B>0$, we have
\eqb\label{eqn:bound_faces}
\mean{\Sigma_{k+1}} = O((\nn r)^{k+2}r^{(d-1)(k+1)}) = O(\nn r) = o(\nn).
\eqe

Therefore, we conclude that
\eqb\label{eqn:additive_l}
\lim_{\nn\to\infty}\frac1\nn \abs{\mean{\Pi_{k,\nn}(\alpha;(0,\lambda]}-\sum_i \meanx{\Pi_{k,\nn}^{(i)}(\alpha;(0,\lambda])}} = 0.
\eqe

Next, note that similarly to the proof of Lemma \ref{lem:pers_assum_linear}, we can choose $\lambda$ large enough so that
\[
\mean{\Pi_{k,\nn}(\alpha;(\lambda,\lmax])}  \le \meanx{F_{k+1,\nn}(\lambda,\lmax)} \approx \nn(F^*_{k+1}-F^\diamond_{k+1}(\lambda))\le \eps\nn,
\]
and
\[
\meanx{\Pi_{k,\nn}^{(i)}(\alpha;(\lambda,\lmax])}  \le \meanx{F_{k+1,\nn}^{(i)}(\lambda,\lmax)} \approx c_i\nn(F^*_{k+1}-F^\diamond_{k+1}(\lambda; M\ind_{Q_i}))\le c_i\eps\nn,
\]
Therefore, for every $\eps>0$ we can take $\lambda$ to be sufficiently large, and have
\[
\lim_{\nn\to\infty}\frac1\nn\abs{\mean{\Pi_{k,\nn}(\alpha;(\lambda,\lmax]}-\sum_i \meanx{\Pi_{k,\nn}^{(i)}(\alpha;(\lambda,\lmax])}} \le 2\eps.
\]
This, combined with \eqref{eqn:additive_l}, concludes the proof.
\end{proof}

\begin{proof}[Proof of Lemma \ref{lem:pers_assum_stable}]
Let $f_1,f_2,\ldots,$ be a sequence of good functions, such that 
\[
\splitb
    \lim_{i\to\infty} f_i(x) &= f(x),\\
    f_i(x) &\le f_j(x),\quad  i<j.
\splite
\]
Define
\[
\delta_i = 1-\int_{\R^d} f_i(x)dx,
\]
so that $\lim_{i\to\infty}\delta_i=0$.
Note that for every $i$ we can write
\[
\cP_\nn = \cP_{\nn,i} \cup \Delta_{\nn,i},
\]
where $\cP_{\nn,i}$ and $\Delta_{\nn,i}$ are independent, and
\[
\splitb
\cP_{\nn,i} &\sim \poisp{(1-\delta_i)\nn}{f_i/(1-\delta_i)}\\
\Delta_{\nn,i} &\sim \poisp{\delta_i \nn}{(f-f_i)/\delta_i}.
\splite
\]
Fix $\alpha>0$ and let $\Pi^{(i)}_{k,\nn}(\alpha)$ be the persistence measure generated by $\cP_{\nn,i}$. Similarly to the previous proofs, we  split the calculations into separate radii ranges. In this case, from Lemma \ref{lem:addone_simplex} and Remark \ref{rem:addone}, we have
\[
\abs{\Pi_{k,\nn}(\alpha; (0,\lambda])) - \Pi^{(i)}_{k,\nn}(\alpha;(0,
\lambda]))} \le \Sigma_{k+1}(\lambda)- \Sigma_{k+1}^{(i)}(\lambda),
\]
where $\Sigma_{k+1}$ 
 and $\Sigma_{k+1}^{(i)}(\lambda)$ count the number of $(k+1)$-simplices generated by $\cP_\nn$ 
 and $\cP_{\nn,i}$ respectively, with radius $\rho$ such that $\nn\rho^d\le \lambda$. Note that
\[
\meanx{\Sigma_{k+1}(\lambda)} = \frac{\nn^{k+2}}{(k+2)!} \int_{(\R^d)^{k+2}}f(\bx)h_r(\bx) d\bx,
\]
where $h_r(\bx)$ is an indicator for $\bx=(x_1,\ldots,x_{k+2})$ spanning a simplex at radius $\rho\le r$.
Taking a change of variables $x_1 \to x$, $x_j \to x+r y_j$, we can show that
\[
\lim_{\nn\to\infty}\frac1\nn \mean{\Sigma_{k+1}(\lambda)} = \frac{\bar h}{(k+2)!} \lambda^{k+1} \int_{\R^d} f^{k+2}(x)dx,
\]
where $\bar h = \int_{(\R^d)^{k+1}} h_1(0,\by)d\by$.
Similarly,
\[
\lim_{\nn\to\infty}\frac1\nn \meanx{\Sigma^{(i)}_{k+1}(\lambda)} = \frac{\bar h}{(k+2)!} \lambda^{k+1} \int_{\R^d} f_i^{k+2}(x)dx.
\]
Using the fact that
\[
f^{k+2}(x)-f_i^{k+2}(x) \le (k+2)\fmax^{k+1}(f(x)-f_i(x)),
\]
we  have,
\eqb\label{eqn:bound_kfaces}
\lim_{\nn\to\infty}\frac1\nn\meanx{\Sigma_{k+1}(\lambda)-\Sigma_{k+1}^{(i)}(\lambda)} \le C \delta_i\lambda^{k+1},
\eqe
for some $C>0$.

Next, using the number of critical faces to bound the death times, we have
\[
\splitb
\abs{\Pi_{k,\nn}(\alpha; (\lambda,\lmax])) - \Pi^{(i)}_{k,\nn}(\alpha;(\lambda,\lmax]))} \le 
F_{k+1,\nn}(\lambda,\lmax) + F_{k+1, \nn}^{(i)}(\lambda,\lmax),
\splite
\]
where $F_{k+1,\nn}$ and $F_{k+1,\nn}^{(i)}$ count the critical faces  in the given range, generated by $\cP_\nn$ and $\cP_{\nn,i}$ respectively.
From Corollary  \ref{cor:crit_lambda_gen}, we have
\[
\lim_{\nn\to\infty}\frac1\nn\meanx{F_{k+1,\nn}^{(i)}(\lambda)} = F_{k+1}^\diamond(\lambda; f_i),
\]
and
\[
\lim_{\lambda\to\infty} F_{k+1}^
\diamond(\lambda; f_i) = (1-\delta_i) F_{k+1}^*.
\]
Therefore, fixing $i_0$, we can find $\lambda_0>0$ such that
\eqb\label{eqn:lb}
F_{k+1}^\diamond(\lambda_0; f_{i_0}) \ge (1-2\delta_{i_0})F_{k+1}^*.
\eqe
Note that the expression for $F_{k+1}^\diamond(\lambda;f)$ in \eqref{eqn:C_lambda_cech} is increasing in $f$, in the sense that
\[
F_{k+1}^\diamond(\lambda;f_{i_0}) \le F_{k+1}^\diamond(\lambda;f_i),\quad  i>i_0,
\]
and therefore, from \eqref{eqn:lb} for all $i\ge i_0$,  we have
\[
\lim_{\nn\to\infty}\frac1\nn\meanx{F_{k+1,\nn}^{(i)}(\lambda_0,\lmax )} =
(1-\delta_i)F_{k+1}^* - F^\diamond_{k+1}(\lambda_0; f_i) \le (2\delta_{i_0}-\delta_{i})F_{k+1}^*.
\]

Finally, fix $\eps>0$, and take $i_0>0$ such that $\delta_{i_0} \le \frac{\eps}{2F_{k+1}^*}$.
Then for all $i\ge i_0$,
\[
\lim_{\nn\to\infty}\frac1\nn\meanx{F_{k+1,\nn}^{(i)}(\lambda_0,\lmax )} 
\le 2\delta_{i_0} F_k^* < \eps.
\]
Taking $\lambda_0$ sufficiently large, also ensures that
\[
\lim_{\nn\to\infty}\frac1\nn \mean{F_{k+1,\nn}(\lambda_0,\lmax)} < \eps,
\]
and going back to \eqref{eqn:bound_kfaces}, for $i$ large enough  we have $C\delta_i\lambda_0^{k+1} < \eps$. Putting it together, we conclude that
\[
\limsup_{\nn\to\infty} \mean{\abs{\Pi_{k,\nn}(\alpha) - \Pi_{k,\nn}^{(i)}(\alpha)}} \le 3\eps.
\]
This concludes the proof.
\end{proof}

Finally, we can prove Theorem \ref{thm:main_poisson}.

\begin{proof}[Proof of Theorem \ref{thm:main_poisson}]
Set $L_\alpha = [\alpha,\infty)$. From Lemma \ref{lem:limit_good}, we conclude that $\frac1\nn\bar\Pi_{k,\nn}(L_\alpha)$ has a universal limit $\Pi_k^*(L_\alpha)$ for all $\alpha\in \R$. Additionally, as the limiting measure is finite, this is sufficient to conclude that $\frac1\nn\Pi_{k,\nn}$ weakly converges to $\Pi_k^*$.
\end{proof}

\subsection{Persistent cycles -- binomial}

In this section we prove Theorem \ref{thm:main_binom}. 

\begin{proof}[Proof of Theorem \ref{thm:main_binom}]
Let $\cX_n\sim\binmp{n}{f}$, and similarly to the previous section, we use the shorthand $\Pi_{k,n}(\alpha) = \Pi_{k,n}([\alpha,\infty))$, for a fixed value of $\alpha$.
Note that we already verified that assumptions \eqref{eqn:assum_translate}-\eqref{eqn:assum_stable} hold for the $\Pi_{k,n}(\alpha)$ functional. Thus, to apply Theorem \ref{thm:univ_binom}, we need to prove the remaining conditions in \eqref{eqn:assum_poly_1}-\eqref{eqn:assum_poly_2}. The main ideas here are inspired by the methods presented in \cite{penrose_random_2003}[Theorem 3.17] and \cite{yogeshwaran_random_2016}[Theorem 4.5].

Using the number of $k$-faces as an upper bound we have $\Pi_{k,n}(\alpha) \le n^{k+1}$, so \eqref{eqn:assum_poly_1} holds.
Denote $\Pi_{k,m,n}(\alpha)$ the number of $\alpha$-persistent $k$-cycles generated by $\cX_{m}$, and with death-time bounded by $\rmax =\rmax(n)$. We need to find an event $A_n$, such that
\[
|\Pi_{k,n}(\alpha)-\Pi_{k,n\pm \Delta,n}(\alpha)| \ind_{A_n} \le C_4 \Delta n^b.
\]
We  first assume that $f$ is of Type I or II, so that $\supp(f)$ is compact. In this case, we can select a collection of boxes $Q_{n,1},\ldots,Q_{n,M}$ whose side-length is $\rmax$, that cover $\supp(f)$, and such that $M \le C\rmax^{-d}$ for some $C>0$.
Defining $N_j = |\cX_{2n}\cap Q_{n,j}|$, then $N_j$ is stochastically dominated by $\mathrm{Binom}(2n,\fmax\rmax^d)$. Therefore, using \cite{penrose_random_2003}[Lemma 1.1], we have for $n$ large enough, for all $j\le M$,
\[
\probx{N_j > n^\delta} \le e^{-n^\delta},
\]
for any $\delta$ such that $n\rmax^d \ll n^\delta$.
Defining  the event
\[
A_{n} := \set{\max_{1\le j \le M} N_j \le n^\delta},
\]
then
\[
\prob{A_{n}^c} \le M e^{-n^\delta} \sim \rmax^{-d} e^{-n^\delta}.
\]
Since $\rmax^d \gg 1/n$, for any $a<\delta$ and large enough $n$,
\[
\prob{A_{n}^c} \le e^{-n^a},
\]
so that \eqref{eqn:assum_A_n} holds. 
Next, define
\[
\Delta_{n,m} = \begin{cases}
    \cX_n \bs \cX_m & n > m,\\
    \cX_m \bs \cX_n & n < m,
\end{cases}
\]
and note that from Lemma \ref{lem:addone_simplex} and Remark \ref{rem:addone},
\[
|\Pi_{k,n}(\alpha)-\Pi_{k,m,n}(\alpha)| \le \Sigma_{k+1},
\]
where $\Sigma_{k+1}$ is the number of $(k+1)$-simplices that include points from $\Delta_{n,m}$, and with radius less than $\rmax$. If $A_{n}$ holds, then all points in $\Delta_{n,m}$ have at most $Cn^\delta$ neighbors, and therefore 
\[
|\Pi_{k,n}(\alpha)-\Pi_{k,n,m}(\alpha)| \le C_2 |n-m| n^{(k+1)\delta},
\]
for some $C_2>0$. Taking $b = (k+1)\delta$, for sufficiently small $\delta$ we have $b<1/4$, and \eqref{eqn:assum_poly_2} holds as well, completing the proof.

Next, suppose that $f$ is of Type III. Fix $\delta$, and take $R=R(n)>0$  such that
\eqb\label{eqn:small_prob}
f(B_R^c) \ll n^{\delta-1},
\eqe
where $B_R$ is a ball of radius $R$ around the origin.
Let $N_R := |\cX_{2n}\cap B_R^c|$, then $N_R\sim\mathrm{Binom}(2n, f(B_R^c))$, and using a similar bound as before, since  $ nf(B_R^c) \ll n^\delta$, we have
\[
    \probx{N_R > n^\delta} \le e^{-n^\delta}.
\]
In addition, we can cover $B_R$ with $M\sim R^d/\rmax^d$ boxes of side-length $\rmax$, 
and we can bound the number of points in each box as before. Therefore, defining
\[
A_{n} := \set{\max(N_1,\ldots, N_M, N_R) \le n^\delta},
\]
then
\[
\prob{A_n^c} \le (M+1)e^{-n^\delta}.
\]

From (III.1), we have that
\[
f(B_R^c) \le C\int_R^\infty \tau^{d-1-A_1}d\tau = C R^{d-A_1}.
\]
Therefore, if we take $R = n^{\frac{\delta-1}{2(d-A_1)}}$, we have that \eqref{eqn:small_prob} holds. Finally, as $M \sim (R/\rmax)^d$, and $n\rmax^d \to\infty$, we have
\[
\prob{A_n^c} \le CnR^d e^{-n^\delta} < e^{-n^a},
\]
for some $a<\delta$. The rest of the proof similar to Types I and II.
\end{proof}

\begin{rem}
    Note that the proof above imposed the condition that $\lmax \ll n^\delta$, while $\delta < \frac1{4(k+1)}$. Thus, for example, if we require $\lmax = n^{\frac1{4(k+2)}}$ the proof would apply.
\end{rem}

\subsection{Critical faces}\label{sec:crit_faces_proofs}

The remaining part to complete is proving 
 the statements for the asymptotics of the number of critical faces presented in Section \ref{sec:crit_faces}.

\subsubsection*{Critical faces in the thermodynamic limit}

\begin{proof}[Proof of Proposition \ref{prop:crit_lambda} -- \v Cech]
Using \eqref{eqn:crit_cech} and Mecke's formula (see \cite{penrose_random_2003}[Theorem 1.6]), we have
\[
\mean{F_{k,\nn}(\lambda)} = \frac{\nn^{k+1}}{(k+1)!} \int_{(\R^d)^{k+1}} f(\bx) h(\bx) \ind\{\rho(\bx)\le r\}e^{-\nn(f(B(\bx)))}d\bx,
\]
where $r = (\lambda/\nn)^{1/d}$, and $h(\bx)$,$\rho(\bx)$, and $B(\bx)$ are defined in Section \ref{sec:crit_faces}, and $f(\bx) = f(x_1)\cdots f(x_{k+1})$.
Next, applying Blaschke-Petkantschin  change of variables \eqref{eqn:bp}, we have
\eqb\label{eqn:crit_td}
 \mean{F_{k,\nn}(\lambda)} = \frac{\nn^{k+1}}{(k+1)!} (k!)^{d-k+1}\int_{\R^d}\int_0^{r} \rho^{dk-1}e^{-\nn f(B_\rho(c))} p(c,\rho) d\rho dc,
\eqe
where
\eqb\label{eqn:pcr}
p(c,\rho):= \int_{\Gamma(d,k)} \int_{(\S^{k-1})^{k+1}} f(c+\rho(\gamma\circ\bth))h(\bth)\vsimp^{d-k+1}(\bth)d\bth d\gamma,
\eqe
$\Gamma(d,k)$ is the $k$-dimensional Grassmannian in $\R^d$, $\S^{k-1}$ is a unit $(k-1)$-dimensional sphere, and $\vsimp(\bth)$ is the volume of the simplex spanned by $\bth$.
Taking the change of variables $\tau = \nn\rho^d$, and taking $\lambda = \nn r^d$, yields
\eqb\label{eqn:F_k_DCT_cech}
\mean{F_{k,\nn}(\lambda)} = \frac{\nn}{d(k+1)!} (k!)^{d-k+1} \int_{\R^d}\int_0^{\lambda} \tau^{k-1}e^{-\nn f(B_{\rho(\tau)}(c))} p(c,\rho(\tau)) d\tau dc,
\eqe
where $\rho(\tau) = (\tau/\nn)^{1/d}$.
Note that the integrand is bounded, and therefore we can apply the dominated convergence theorem (DCT). Taking the limit, we have
\[
\lim_{\nn\to\infty} p(c,(\tau/\nn)^{1/d}) = f^{k+1}(c) V_{d,k} \Gamma_{d,k},
\]
where $V_{d,k} = \int_{(\S^{k-1})^{k+1}} h(\bth)\vsimp^{d-k+1}(\bth)d\bth$, and $\Gamma_{d,k}$ is the volume of the Grassmannian.
Also, using the Lebesgue differentiation theorem, we have
\[
\lim_{\nn\to\infty} \nn f(B_{\rho(\tau)}(c)) = \omega_d \tau f(c)  ,
\]
where $\omega_d$ is the volume of a $d$-dimensional unit ball.
Thus, we have
\[
\splitb
\lim_{\nn\to\infty} \frac1\nn \mean{F_{k,\nn}(\lambda)} &= \frac{V_{d,k}\Gamma_{d,k}(k!)^{d-k+1}}{d(k+1)!}\int_{\R^d}f^{k+1}(c)\int_0^\lambda \tau^{k-1} e^{-\omega_d \tau f(c)}d\tau dc\\
&= \frac{V_{d,k}\Gamma_{d,k}(k!)^{d-k+1}}{d\omega_d^k(k+1)!}\int_{\R^d}f(c)\int_0^{\omega_d f(c)\lambda} t^{k-1} e^{-t}dt dc.
\splite
\]
We can write this limit as
\eqb\label{eqn:C_lambda_cech}
F_k^\diamond(\lambda;f) := \frac{V_{d,k}\Gamma_{d,k}(k!)^{d-k+1}}{d\omega_d^k(k+1)!} \int_{\R^d} f(c)\gamma^{_\cC}_k(\omega_d f(c)\lambda)dc,
\eqe
where $\gamma^{_\cC}_k$ is the lower incomplete gamma function. Noting that $\lim_{\lambda\to\infty} \gamma^{_\cC}_k(\omega_d f(c)\lambda) = \Gamma(k) = (k-1)!$, we have 
\eqb\label{eqn:C_infty_cech}
F_k^* = \lim_{\lambda\to\infty} F_k^\diamond(\lambda;f) = \frac{V_{d,k}\Gamma_{d,k}(k!)^{d-k+1}}{d\omega_d^kk(k+1)}, 
\eqe
which is independent of $f$. This completes the proof. 
\end{proof}

\begin{proof}[Proof of Proposition \ref{prop:crit_lambda} --  Vietoris-Rips]
Using \eqref{eqn:crit_rips} and Mecke's formula, for $k\ge2$  we have
\[
    \mean{F_{k,\nn}(\lambda)} =
\frac{\nn^2}{2}\int_{(\R^d)^2} f(x_1)f(x_2)\mean{\beta_{k-2}(x_1,x_2)}\ind\set{|x_1-x_2|\le r} dx_1 dx_2,
\]
where 
\[
\beta_{k-2}(x_1,x_2) := \beta_{k-2}(\lk((x_1,x_2);\cP_\nn\cup\set{x_1,x_2})).
\]
Taking the change of variables $x_1\to x$, $x_2 \to x+\rho\theta$, we have
\eqb\label{eqn:F_k_rips_polar}
\mean{F_{k,\nn}(\lambda)} = \frac{\nn^2}{2}\int_{\R^d}\int_0^{r} \int_{\S^{d-1}} f(x)f(x+\rho\theta) \rho^{d-1} \mean{\beta_{k-2}(x,x+\rho\theta)} d\theta d\rho dx.
\eqe
Taking $\tau = \nn\rho^d$, we have
\eqb\label{eqn:C_k_rips}
\mean{F_{k,\nn}(\lambda)} = \frac{\nn}{2d} \int_{\R^d}\int_0^{\lambda}\int_{\S^{d-1}}f(x)f(x+\rho(\tau)\theta)\meanx{\beta_{k-2}(x,x+\rho(\tau)\theta)}d\theta d\tau dx,
\eqe
with $\rho(\tau) = (\tau/\nn)^{1/d}$.
Denote  $I(x_1,x_2) = B_{|x_1-x_2|}(x_1) \cap B_{|x_1-x_2|}(x_2)$. By bounding the number of $(k-1)$-faces, we have
\[
\beta_{k-2}(x_1,x_2) \le |\cP_\nn\cap I(x_1, x_2)|^{k-1}.
\]
Since the number of points is a Poisson random variable, we have
\[
\mean{\beta_{k-2}(x_1,x_2)} \le C (\nn f(I(x_1,x_2)))^{k-1},
\]
for some $C>0$. Note that
\[
f(I(x_1,x_2)) \le \fmax \vol(I(x_1,x_2)) = \fmax |x_1-x_2|^{d} \kappa_d,
\]
where $\kappa_d = \vol(I_0)$, and $I_0 = B_1({\mathbf 0})\cap B_1(\mathbf{1})$, with $\mathbf{0}$ being the origin in $\R^d$, and $\mathbf{1} = (1,0,\ldots,0)\in \R^d$.
Thus,
\[
\meanx{\beta_{k-2}(x,x+\rho(\tau)\theta)} \le C( \tau\fmax \kappa_d)^{k-1},
\]
and we can apply the DCT to \eqref{eqn:C_k_rips}. To find the limit of the integrand, we note that
\eqb\label{eqn:betti_sum}
\mean{\beta_{k-2}(x_1,x_2)} = e^{-\nn f(I(x_1,x_2))}\sum_{m=2(k-1)}^\infty \frac1{m!}(\nn f(I(x_1,x_2)))^m \bar\beta_{k-2}(x_1,x_2;m), 
\eqe
where 
\[
\bar \beta_{k-2}(x_1,x_2;m) := \cmean{\beta_{k-2}(x_1,x_2)}{|\cP_\nn \cap I(x_1,x_2)|=m}.
\]
For $x_1=x$ and $x_2 = x+\rho(\tau)\theta$, we can bound the term inside the summation by
\[
\frac1{m!} (\tau\fmax \kappa_d)^m m^{k-1},
\]
which is summable, and therefore we can use the DCT for the sum as well.
Note that $\lk((x_1,x_2);\cP_\nn\cup\set{x_1,x_2}) = \cR_{|x_1-x_2|}(\cP_\nn\cap I(x_1,x_2))$ (i.e., the link itself is a Vietoris-Rips complex), and therefore,
\[
\bar\beta_{k-2}(x_1,x_2;m) = \int_{(I(x_1,x_2))^m} \frac{f(\by)}{(f(I(x_1,x_2)))^m}\beta_{k-2}(\cR_{|x_1-x_2|}(\by)) d\by.
\]
Next, we take  $x_1=x, x_2 = x+\rho\theta$, and the change of variables $y_i \to x+\rho\cdot( \theta\circ z_i)$, where $z_i\in I_0$ and $\theta\circ z_i$ is the rotation in the direction of $\theta\in \S^{d-1}$. Then,
\[
\bar\beta_{k-2}(x,x+\rho\theta ;m) = (f(I(x,x+\rho\theta)))^{-m}\rho^{dm}\int_{I_0^m} f(x+\rho\cdot(\theta\circ\bz))\beta_{k-2}(\cR_{1}(\bz)) d\bz,
\]
where we used the fact that $\beta_{k-2}(\cR_\rho(x+\rho(\theta\circ\bz))) = \beta_{k-2}(\cR_1(\bz))$, by translation and rotation invariance and scaling. Putting this back into \eqref{eqn:betti_sum} yields,
\[
\meanx{\beta_{k-2}(x,x+\rho(\tau)\theta)} = e^{-\nn f(I(x,x+\rho(\tau)\theta))}\sum_{m=2(k-1)}^\infty \frac{\tau^m}{m!}\int_{I_0^m} f(x+\rho(\tau)(\theta\circ \bz))\beta_{k-2}(\cR_1(\bz))\bz,
\]
and when $\nn\to\infty$, since $\rho(\tau)\to 0$, we have
\[
\lim_{\nn\to\infty}\meanx{\beta_{k-2}(x,x+\rho(\tau)\theta)} = e^{-\tau f(x) \kappa_d}\sum_{m=2(k-1)}^\infty \frac{(\tau f(x) \kappa_d)^m}{m!}\beta_{k-2}^*(m),
\]
where
\[ \beta_{k-2}^*(m) := \int_{I_0^m}\kappa_d^{-m}\beta_{k-2}(\cR_1(\bz))\bz.
\]
is the expected reduced Betti number for $m$ points uniformly distributed in $I_0$. Putting back into \eqref{eqn:C_k_rips}, we have
\[
\lim_{\nn\to\infty}\frac1\nn \mean{F_{k,\nn}(\lambda)} = \frac{\omega_d}{2} \int_{\R^d} \int_0^{\lambda}f^2(x) e^{-\tau f(x)\kappa_d}\param{\sum_{m=2(k-1)}^\infty \frac{(\tau f(x) \kappa_d)^m }{m!}\beta_{k-2}^*(m)}d\tau dx.
\]
For a fixed $x$, we can take a change of variables $t=\tau f(x) \kappa_d$, and then,
\[
\lim_{\nn\to\infty}\frac1\nn \mean{F_{k,\nn}(\lambda)} = \frac{\omega_d}{2\kappa_d} \int_{\R^d} f(x) \int_0^{\lambda f(x) \kappa_d} e^{-t}\param{\sum_{m=2(k-1)}^\infty \frac{t^m }{m!}\beta_{k-2}^*(m)}dt dx.
\]
Defining
\[
\gamma^{_\cR}_k(z) := \int_0^{z} e^{-t}\param{\sum_{m=2(k-1)}^\infty\frac{t^m }{m!}\beta_{k-2}^*(m)}dt, 
\]
then
\eqb\label{eqn:C_lambda_rips}
F_k^\diamond(\lambda;f) :=\frac{\omega_d}{2\kappa_d} \int_{\R^d} f(x) \gamma_k^{_\cR}(\lambda f(x) \kappa_d)dx.
\eqe
Note that
\[
\gamma_k^{_\cR}(\infty) = \sum_{m=2(k-1)}^\infty\beta^*_{k-2}(m).
\]
Therefore,
\eqb\label{eqn:C_infty_rips}
F_k^* = \lim_{\lambda\to\infty}F_k^\diamond(\lambda;f) =\frac{\omega_d}{2\kappa_d}\sum_{m=2(k-1)}^\infty\beta^*_{k-2}(m),
\eqe
which is indeed independent of $f$, completing the proof for $k\ge 2$.

The case $k=1$ is simpler. Recall from Lemma \ref{lem:crit_rips} that here we only need to count isolated edges, i.e., such that $\lk(e;\cP_\nn) = \emptyset$. Therefore, 
\[
\splitb \mean{F_{1,\nn}(\lambda)} &=
\frac{\nn^2}{2}\int_{(\R^d)^2} f(x_1)f(x_2)\prob{\cP_\nn\cap I(x_1,x_2) = \emptyset}\ind\set{|x_1-x_2|\le r} dx_1 dx_2,\\
&=\frac{\nn^2}{2}\int_{(\R^d)^2} f(x_1)f(x_2)e^{-\nn f(I(x_1,x_2))}\ind\set{|x_1-x_2|\le r} dx_1 dx_2.
\splite
\]
Similar limits as above, yield
\[
\splitb
\lim_{\nn\to\infty}\frac1\nn\mean{F_{1,\nn}(\lambda)} &= F_1^\diamond(\lambda;f) := \frac{\omega_d}{2\kappa_d} \int_{\R^d} f(x)(1-e^{-\lambda f(x) \kappa_d})dx,\\
\lim_{\lambda\to\infty} F_{1}^\diamond(\lambda;f) &= F_1^*(\lambda) := \frac{\omega_d}{2\kappa_d}.
\splite
\]
\end{proof}

\begin{proof}[Proof of Corollary \ref{cor:crit_lambda_gen}]
Recall that here $\cP_\nn\sim\poisp{c\nn}{\tilde f}$, where $\tilde f = f/c$. In this case, for the \v Cech complex, \eqref{eqn:crit_td} should be replaced with
\eqb\label{eqn:crit_td_gen}
\mean{F_{k,\nn}(\lambda)} = \frac{(c\nn)^{k+1}}{(k+1)!}(k!)^{d-k+1} \int_{\R^d}\int_0^{r} \rho^{dk-1}e^{-c\nn \tilde f(B_\rho(c))} \tilde p(c,\rho) d\rho dc,
\eqe
where
\[
\splitb
\tilde p(c,\rho) &:= \int_{\Gamma(d,k)} \int_{(\S^{k-1})^{k+1}} \tilde f(c+\rho(\gamma\circ\bth))h(\bth)\vsimp^{d-k+1}(\bth)d\bth d\gamma\\
&= c^{-(k+1)}
\int_{\Gamma(d,k)} \int_{(\S^{k-1})^{k+1}} f(c+\rho(\gamma\circ\bth))h(\bth)\vsimp^{d-k+1}(\bth)d\bth d\gamma,
\splite
\]
and 
\[
c\tilde f(B_\rho(c)) = f(B_\rho(c)).
\]
Thus, \eqref{eqn:crit_td_gen} is equal to the right hand side of \eqref{eqn:crit_td}. So we can continue with the same steps as the proof of Proposition \ref{prop:crit_lambda}, to have
\[
\lim_{\nn\to\infty}\frac1\nn\mean{F_{k,\nn}(\lambda)} = F_k^\diamond(\lambda).
\]
Note that in this case, taking the limit in \eqref{eqn:C_lambda_cech}, yields
\[
\lim_{\lambda\to\infty} F_k^\diamond(\lambda) = \frac{V_{d,k}\Gamma_{d,k}(k!)^{d-k+1}}{d\omega_d^k k(k+1)} \int_{\R^d} f(c)dc = c F_k^*.
\]
Similar arguments apply for the Vietoris-Rips case.
\end{proof}

\subsubsection*{The total number of critical faces}


\begin{proof}[Proof of Lemma \ref{lem:good} -- \v Cech]

Take $r\to 0$ and $\Lambda = \nn r^d\to\infty$.
Similarly to \eqref{eqn:crit_td}, we have
\eqb\label{eqn:ckr_vol}
 \splitb
\mean{F_{k,\nn}(\Lambda)} 
&=\frac{\nn^{k+1}}{(k+1)!} (k!)^{d-k+1}\int_{\R^d}\int_0^{r} \rho^{dk-1} e^{-\nn f(B_\rho(c))} p(c,\rho)d\rho dc\\
&= \frac{\nn^{k+1}}{d\omega_d^k(k+1)!} (k!)^{d-k+1} \int_{\R^d}\int_0^r \Delta(c,\rho) f(\S_\rho(c))f^{k-1}(B_\rho(c)) e^{-\nn f(B_\rho(c))}d\rho dc,
\splite
\eqe
where
\eqb\label{eqn:delta}
\Delta(c,\rho) := 
\frac{d\omega_d^k\rho^{dk-1}p(c,\rho)}{f(S_\rho(c))f^{k-1}(B_\rho(c))} = \frac{p(c,\rho)}{\bar f(\S_\rho(c))\bar f^{k-1}(B_\rho(c))},
\eqe
$\S_\rho(c)$ is the $(d-1)$-sphere centered at $c$, and $\bar f$ is the average value. 
Note that by $f(\S_\rho(c))$ we refer to the integral with respect to the $(d-1)$-dimensional surface measure of the sphere.

Fix $c$, and denote $f_c(\rho) = f(B_\rho(c))$. Then $f_c(\rho)$ is increasing (in $\rho$), and $\frac{d}{d\rho}f_c(\rho) = f(S_\rho(c))$. We will thus take the change of variables $t = \nn f_c(\rho)$, leading to
\eqb\label{eqn:mean_trans}
\mean{F_{k,\nn}(\Lambda)} = \frac{\nn}{d\omega_d^k(k+1)!} (k!)^{d-k+1}\int_{\R^d}\int_0^\infty \Delta(c,\rho(t)) \indf{t\le \nn f_c(r)} t^{k-1}e^{-t} d\rho dc,
\eqe
where $\rho(t) = f_c^{-1}(t/\nn)$. Note that for any fixed $c\in \supp(f)$ and $t<\nn$, we have
\[
\lim_{\nn\to \infty}\rho(t) = 0,
\]
and
\[
\lim_{\rho\to 0}\Delta(c,\rho) = {\Gamma_{d,k}V_{d,k}f(c)}.
\]
Therefore, if we are allowed to take the point-wise limit of the integrand in \eqref{eqn:mean_trans} we get
\[
\lim_{\nn\to\infty} \frac1\nn\mean{F_{k,\nn}(\Lambda)} = F_k^* 
\]
as in \eqref{eqn:C_infty_cech}.
Our goal next will be to justify taking the limit for the integrand in \eqref{eqn:mean_trans}, by bounding $\Delta(c,\rho)$, and applying the DCT. The approach for the bounds will be different between the three types presented in Section \ref{sec:prob_set}.

\noindent{\bf Type I:} Recall that here the support $\cS = \supp(f)$ is compact and $\fmin = \inf_{\cS} f >0$.
Let $\partial \cS$ be the boundary of $\cS$, and define
\[
\cS^{(r)} := \set{x\in \cS : \inf_{z\in \partial \cS}|x-z|>r},
\]
i.e., all points in $\cS$ that are at least $r$ away from the boundary.
Next, we split $F_{k,\nn}(\Lambda)$ into $F_{k,\nn}^{(1)}(\Lambda)+F_{k,\nn}^{(2)}(\Lambda)$, where $F_{k,\nn}^{(1)}(\Lambda)$ counts critical $k$-faces whose center $c$ is in $\cS^{(r)}$, and $F_{k,\nn}^{(1)}(\Lambda)$ counts the rest.
If $c\in \cS^{(r)}$, and $\rho \le r$, then $B_\rho(c) \subset \cS$, and thus
\[
\bar f(B_\rho(c))\ge \fmin \quad\text{and}\quad \bar f(S_\rho(c))\ge \fmin.
\]
Also, note that from \eqref{eqn:pcr}
\[
\splitb
p(c,\rho)&\le V_{\max}\int_{\Gamma(d,k)} \param{\int_{\S^{k-1}} f(c+\rho\gamma\theta)d\theta}^{k+1} d\gamma\\
&\le V_{\max}\Gamma_{d,k}(k\omega_k\fmax)^{k+1},
\splite
\]
where $V_{\max} = \sup_{(\S^{k-1})^{k+1}}\vsimp^{d-k+1}(\bth)$, and $\fmax = \sup_{\cS}f(x)< \infty$. Therefore, from \eqref{eqn:delta}, we have
\[
\Delta(c,\rho) \le C = \frac{V_{\max}\Gamma_{d,k}(k\omega_k\fmax)^{k+1}}{\fmin^k}.
\]
Similarly to 
\eqref{eqn:mean_trans}, we have
\[
\meanx{F_{k,\nn}^{(1)}(\Lambda)} = \frac{\nn}{d\omega_d^k(k+1)!}(k!)^{d-k+1} \int_{\R^d}\int_0^\infty \Delta(c,\rho(t)) \indf{t\le \nn f_c(r)} \indfx{c\in \cS^{(r)}} t^{k-1}e^{-t} d\rho dc,
\]
where the the integrand is bounded and we can apply the DCT. Note that for every $c\in \cS\bs\partial \cS$ we have $\ind\{c\in \cS^{(r)}\}\to 1$. Therefore,
\[
\lim_{\nn\to\infty}\frac1\nn\meanx{F_{k,\nn}^{(1)}(\Lambda)} = F_k^*.
\]
To bound $F_{k,\nn}^{(2)}$, we note that if $c(\bx)\in \cS^{(r)}$ then all the points in $\bx$ must lie in $\cS\bs \cS^{(2r)}$. Denote by $N_{k,\nn}$ the number of $k$-faces in $\cC_r(\cP_{\nn})$ that are contained in $\cS\bs \cS^{(2r)}$. Then,
\[
\splitb
\meanx{F_{k,\nn}^{(2)}(\Lambda)} \le \mean{N_{k,\nn}} &\le \frac{\nn^{k+1}}{(k+1)!} \int_{(\cS\bs \cS^{(2r)})^{k+1}}f(\bx)\ind\set{|x_i-x_1| \le 2r,\ 2\le i \le k+1}d\bx\\ &\le C \nn^{k+1} r^{dk+1},
\splite
\]
for some $C>0$, where we used the fact that $\vol(\cS\bs \cS^{(2r)}) = O(r)$.
Therefore,
\[
\frac1\nn\meanx{F_{k,\nn}^{(2)}(\Lambda)} \le C r (\nn r^d)^k  = C \param{\frac{\Lambda}{\nn}}^{1/d} \Lambda^k.
\]
Thus, taking $\lmax = \nn^{\frac1{dk+2}}$,
 we have
\[
\frac1\nn\meanx{F_{k,\nn}^{(2)}(\lmax)}  \to 0,
\]
implying that
\[
\lim_{\nn\to\infty}\frac1\nn\mean{F_{k,\nn}(\lmax)}  = F_k^*,
\]
completing the proof for Type I.

\noindent{\bf Type II:} 
Recall that here $\cS$ is compact, while $\fmin = 0$. We defined $\delta_0$ so that
for all $x\in \cS$ at most $\delta_0$ away from $\cS_0 = f^{-1}(0)\cap \cS$, we have $C_1(\delta(x))^q \le f(x) \le C_2(\delta(x))^q$.

Consider first, the case where $c\in \cS^{(2r)}\bs \cS^{(\delta_0-r)} $. Defining
\eqb\label{eqn:fminmax}
\fmin(c,\rho) = \inf_{B_\rho(c)} f,\quad \fmax(c,\rho) = \sup_{B_\rho(c)}f,
\eqe
we have
\[
\fmax(c,\rho) \le C_2(\delta(c)+r)^\alpha, \qquad \fmin \ge C_1(\delta(c)-r)^\alpha,
\]
and therefore,
\[
\frac{\fmax(c,\rho)}{\fmin(c,\rho)}\le \frac{C_2}{C_1}\param{\frac{\delta(c)+r}{\delta(c)-r}}^{\alpha} =  \frac{C_2}{C_1}\param{\frac{1+r/\delta(c)}{1-r/\delta(c)}}^{\alpha}
\]
Since $\delta(c)\ge 2r$, the last term is bounded by a constant $C>0$.

Next, consider $c\in \cS^{(\delta_0-r)}$. In this case, we have $B_\rho(c)\subset \cS^{(\delta_0-2r)}$. For $r$ small enough, $\delta_0-2r > \delta_0/2$, where $\cS^{(\delta_0/2)}$ is compact, and $\tilde{f}_{\min} := \inf_{\cS^{(\delta_0/2)}} f >0$. Therefore,
\[
	\bar f(B_\rho(c))\ge \tilde{f}_{\min}\quad\text{and}\quad \bar f(\S_\rho(c))\ge \tilde{f}_{\min}
\]
and
\[
	\Delta(c,\rho) \le \frac{V_{\max} \Gamma_{d,k} (k\omega_k \fmax)^{k+1}}{\tilde{f}_{\min}^k}.
\]
To conclude, we have that for all $c\in \cS^{(2r)}$, $\Delta(c,\rho)$ is upper bounded by a constant $C$.
We can therefore continue the same way we did for Type I to complete the proof.

\noindent{\bf Type III:}
Here, $\cS = \R^d$, and $f(x) = c e^{-\Psir(x)\Psis(x/|x|)}$.
First, we write
\[
F_{k,\nn}(\Lambda) = F_{k,\nn}^{(1)}(\Lambda) + F_{k,\nn}^{(2)}(\Lambda) + F_{k,\nn}^{(3)}(\Lambda),
\]
with 
\[
\splitb
F_{k,\nn}^{(1)}(\Lambda) &:= \# \text{critical $k$-faces with } \nn\rho^d \le \Lambda \text{ and } |c| \le R_0,
\\
F_{k,\nn}^{(2)}(\Lambda) &:= \# \text{critical $k$-faces with } \nn\rho^d \le \Lambda \text{ and } R_0 < |c| \le R ,\\
F_{k,\nn}^{(3)}(\Lambda) &:= \# \text{critical $k$-faces with } \nn\rho^d \le \Lambda \text{ and } R < |c|,
\splite
\]
where $R_0$ is the radius required for conditions (III.1)-(III.3) to hold, and $R$ will be chosen later.
For $F_{k,\nn}^{(1)}$ note that $B_{R_0}(\origin)$ is compact, with $\inf_{B_{R_0}(\origin)}f>0$. Therefore,  similarly to Type I, we have,
\eqb\label{eqn:FK1}
\lim_{\nn\to\infty}\frac1\nn\meanx{F_{k,\nn}^{(1)}(\Lambda)} = f(B_{R_0}(\origin))F_k^*.
\eqe
Next, we skip to $F_{k,\nn}^{(3)}$.  Denote by $N_{k,\nn}$ the number of $k$-faces in $\cC_r(\cP_\nn)$, with at least one vertex lying in $\R^d\bs B_R(\origin)$. Then
\[
\splitb 
\meanx{F_{k,\nn}^{(3)}(\Lambda)} \le \mean{N_{k,\nn}} &\le (k+1)\frac{\nn^{k+1}}{k!} \int_{(\R^d\bs B_R(\origin))\times(\R^d)^k} f(\bx) \ind\set{|x_i-x_1| \le 2r,  2\le i \le k+1}d\bx\\
&\le C\nn^{k+1}r^{dk}\int_R^\infty \int_{\S^{d-1}} \tau^{d-1} e^{-\Psir(\tau)\Psis(\theta)}d\theta d\tau\\
&\le C\nn^{k+1}r^{dk}\int_R^\infty \tau^{d-1} e^{-\Psir(\tau)} d\tau,
\splite
\]
since we assume $\Psis\ge 1$.
Take the change of variables $t=\Psir(\tau)$, then
\[
\meanx{F_{k,\nn}^{(3)}(\Lambda)}  \le C\nn^{k+1}r^{dk} \int_{\Psir(R)}^\infty\frac{(\Psir^{-1}(t))^{d-1}}{\Psir'(\Psir^{-1}(t))}e^{-t}dt.
\]
Using assumption (III.2), we have
\[
\meanx{F_{k,\nn}^{(3)}(\Lambda)} \le C\nn^{k+1}r^{dk} \int_{\Psir(R)}^\infty(\Psir^{-1}(t))^{d}e^{-t}dt.
\]
Finally, using assumption (III.1), we have 
$
\Psir^{-1}(t)\le e^{t/A_1} 
$.
Therefore,
\[
\meanx{F_{k,\nn}^{(3)}(\Lambda)}  \le C\nn^{k+1}r^{dk}\int_{\Psir(R)}^\infty e^{-t\param{1-\frac{d}{A_1}}}dt 
= C\nn^{k+1}r^{dk} e^{-\Psir(R)\param{1-\frac{d}{A_1}}}.
\]
Thus, if $A_1 > d$ we can take
\eqb\label{eqn:R}
R = \Psir^{-1}\param{\frac{k+1}{1-d/A_1}\log\nn},
\eqe
and then  $\meanx{F_k^{(3)}(\Lambda)}\to 0$.

The final step is to evaluate $\meanx{F_{k,\nn}^{(2)}(\Lambda)}$.
We go back to \eqref{eqn:ckr_vol}, and try to bound \eqref{eqn:delta}. Recall, that
\[
p(c,\rho) =\int_{\Gamma(d,k)} \int_{(\S^{k-1})^{k+1}} f(c+\rho(\gamma\circ\bth))h(\bth) \vsimp^{d-k+1}(\bth)d\bth d\gamma.
\]
The condition imposed by $h(\bth)$ implies that at least one of the points $z = c+\rho\gamma\theta_i$ satisfies $|z| > |c|$ (otherwise all points are on a single hemisphere, and the simplex is not critical). Recall that here we have $|c| > R_0$, therefore from monotonicity (III.2),
\[
\Psir(|z|) \ge \Psir(|c|),
\]
and from the Lipschitz condition,
\[
|\Psis(z/|z|)-\Psis(c/|c|)| \le K\frac\rho{|c|},
\]
for some $K>0$.
This implies that
\[
\frac{f(z)}{f(c)} \le e^{\Psir(|c|)(\Psis(c/|c|)-\Psis(z/|z|))}\le e^{K\Psir(|c|)\frac{\rho}{c}} \le e^{K\Psir(R)\frac{r}{R_0}}.
\]
Therefore, recalling \eqref{eqn:R} we have that taking any $r=O\param{\frac1{\log\nn}}$, 
\[
p(c,\rho) \le C f(c) \fmax^k(c,\rho),
\]
for some $C>0$.
Thus,
\[
\Delta(c,\rho) \le C f(c)\param{\frac{\fmax(c,\rho)}{\fmin(c,\rho)}}^{k}.
\]

Next, let $z = c+\rho\theta\in B_\rho(c)$. Then, using assumption (III.3), we have
\[
\Psir(|c|)(1-A_3r)\le \Psir(|c|-r)\le  
\Psir(|z|) \le \Psir(|c|+r) \le \Psir(|c|)(1 + A_3 r)
\]
In addition, since $\Psis$ is Lipschitz, we have
\[
\Psis(c/|c|)- K\frac{r}{R_0}\le \Psis(z/|z|) \le \Psis(c/|c|) + K\frac{r}{R_0}.
\]
Therefore,
\[
\frac{\fmax(c,\rho)}{\fmin(c,\rho)} \le e^{2r\Psir(R)\param{\frac{K}{R_0} + A_3\Psis(c/|c|)}} \le e^{Cr\Psir(R)},
\]
since $\Psis$ is bounded. Therefore, taking $r = O\param{\frac{1}{\log\nn}}$, from \eqref{eqn:R} we have that,
\[
\Delta(c,\rho) \le C f(c).
\]
This  implies that we can evaluate $\meanx{F_{k,\nn}^{(2)}(\Lambda)}$ similarly to \eqref{eqn:mean_trans},
\[
\splitb
\meanx{F_{k,\nn}^{(2)}(\Lambda)} &= \frac{\nn}{d\omega_d^k(k+1)!}(k!)^{d-k+1}\\
&\times\int_{\R^d}\int_0^\infty \Delta(c,\rho(t)) \indf{t\le \nn f_c(r)} \ind\set{R_0<|c|\le R}t^{k-1}e^{-t} d\rho dc,
\splite
\]
where we have an integrable bound for the integrand.
Applying the DCT as before, assuming $\rmax = O(1/\log \nn)$ and $\lmax = \nn\rmax^d$, we have
\[
\lim_{\nn\to\infty}\frac1\nn\meanx{F_{k,\nn}^{(2)}(\lmax)} = f(\R^d\bs B_{R_0}(\origin)) F_k^*.
\]
Together with \eqref{eqn:FK1},this 
completes the proof.
\end{proof}

\begin{proof}[Proof of Lemma \ref{lem:good} -- Vietoris-Rips]
Recall \eqref{eqn:F_k_rips_polar}, for $k\ge 2$,
\[
\mean{F_{k,\nn}(\Lambda)} = \frac{\nn^2}{2} \int_{\R^d}\int_0^{r}\int_{\S^{d-1}}\rho^{d-1}f(x)f(x+\rho\theta)\meanx{\beta_{k-2}(x,x+\rho\theta)}d\theta d\rho dx,
\]
where from \eqref{eqn:betti_sum},
\[
\meanx{\beta_{k-2}(x,x+\rho\theta)}= \sum_{m=2(k-1)}^\infty \frac1{m!} {(\nn f(I(x,x+\rho\theta))^m}e^{-\nn f(I(x,x+\rho\theta))}\bar\beta_{k-2}(x,x+\rho\theta;m).
\]
Denote $f_{x,\theta}(\rho) = f(I(x,x+\rho\theta))$,  and $\Delta(x,\rho,\theta) :=
\frac{\rho^{d-1}f(x+\rho\theta)}{f'_{x,\theta}(\rho)}
$, then we can write
\[
\splitb
\mean{F_{k,\nn}(\Lambda)} = \frac{\nn}{2}\int_{\R^d}\int_0^{r}\int_{\S^{d-1}} f(x)\Delta(x,\rho,\theta) \sum_{m=2(k-1)}^\infty &\frac{1}{m!} \bar\beta_{k-2}(x,x+\rho\theta;m)\\
&\times (\nn f'_{x,\theta}(\rho)) (\nn f_{x,\theta}(\rho))^m e^{-\nn f_{x,\theta}(\rho)}d\theta d\rho dx.
\splite
\]

Taking the change of variables $t= \nn f_{x,\theta}(\rho)$ and $dt = \nn f'_{x,\theta}(\rho)d\rho$, we have
\eqb\label{eqn:C_k_simp}
\splitb
\mean{F_{k,\nn}(\Lambda)} = \frac{\nn}{2}\int_{\R^d}\int_0^\infty\int_{\S^{d-1}} &f(x)\Delta(x,\rho(t),\theta) \indf{t\le \nn f_{x,\theta}(r)}\\
 &\times {\sum_{m=2(k-1)}^\infty \frac{1}{m!} \bar\beta_{k-2}(x,x+\rho(t)\theta;m)
 t^m e^{-t}} d\theta dt dx,
 \splite
\eqe
where
\[
\rho(t) := f_{x,\theta}^{-1}(t/\nn).
\]
For all three types of distributions, we saw in the proof for the \v Cech complex that we can upper bound the term $\fmax(x,\rho)/\fmin(x,\rho)$. We will use that to upper bound the integrand here as well.
The rest of the details will be the same as for the \v Cech complex.

First, note that 
\[
f'_{x,\theta}(\rho) = \lim_{\eps\to 0}\frac{f_{x,\theta}(\rho+\eps)-f_{x,\theta}(\rho)}{\eps} =   \lim_{\eps\to 0}\frac1\eps \int_{I(x,x+(\rho+\eps)\theta)\bs I(x,x+\rho\theta)}f(z)dz \ge \fmin(x,\rho) V_{x,\theta}'(\rho),
\]
where $V_{x,\theta}(\rho) = \vol(I(x,x+\rho\theta)) = \rho^d \vol(I_0)$, and therefore, $V_{x,\theta}'(\rho) = d\rho^{d-1} \kappa_d$.
Thus,
\[
f'_{x,\theta}(\rho) \ge d\kappa_d \rho^{d-1}\fmin(x,\rho),
\]
and
\eqb\label{eqn:rips_int_bound_1}
\Delta(x,\rho,\theta)\le \frac1{d\kappa_d}\frac{\fmax(x,\rho)}{\fmin(x,\rho)} \le C,
\eqe
for some $C>0$ (using the bounds we established for the \v Cech case).

Next, Lemma~\ref{lem:forbidden} provides a necessary condition for
$\beta_{k-2}(x_1,x_2)$ to be nonzero. This condition requires that a ``forbidden region" $I_\varobslash(x_1,x_2)\subset I(x_1,x_2)$ may contain no points from $\cP_\nn$.
Therefore,
\[
\bar\beta_{k-2}(x_1,x_2;m) \le \param{1- \frac{f(I_{\varobslash}(x_1,x_2))}{f(I(x_1,x_2))}}^m m^{k-1}. 
\]
Then,
\[
\bar \beta_{k-2}(x,x+\rho\theta;m) \le \param{1-\frac{\fmin(x,\rho)\vol(I_{\varobslash}(x,x+\rho\theta))}{\fmax(x,\rho)\vol(I(x,x+\rho\theta))}}^m m^{k-1} \le \xi^m m^{k-1},
\]
for some constant $0 < \xi < 1$, since the volume of $I_\varobslash$ is proportional to $I$ (see Lemma \ref{lem:forbidden}). In this case,
\eqb\label{eqn:rips_int_bound_2}
\splitb
\sum_{m=2(k-1)}^\infty \frac{1}{m!} \bar\beta_{k-2}(x,x+\rho(t)\theta;m)
 t^m e^{-t} &\le \sum_{m=2(k-1)}^\infty \frac{1}{m!} 
 m^{k-1}(\xi t)^m e^{-t} \\
 &\le C(\xi t)^{k-1} e^{-(1-\xi)t},
\splite
\eqe
for some $C>0$.
Combining \eqref{eqn:rips_int_bound_1} and \eqref{eqn:rips_int_bound_2}, we can apply the same steps as in the \v Cech complex, and have
\[
\lim_{\nn\to\infty}\frac1\nn\mean{F_{k,\nn}(\Lambda)} = F_k^*.
\]
\end{proof}

\begin{proof}[Proof of Lemma \ref{lem:comp_conv}]

For the \v Cech filtration, note that in the case where the support is compact and convex, and the density is bounded away from zero, we have $f(B_\rho(c)) \ge C \rho^d$, 
for some constant $C>0$. Thus, the integrand in 
\eqref{eqn:F_k_DCT_cech} is bounded by an integrable term, even if we take $\lambda=\infty$, and applying the DCT to \eqref{eqn:F_k_DCT_cech} with $\lambda=\infty$ leads to the desired result. Similar arguments apply for the Vietoris-Rips complex, by applying the DCT to \eqref{eqn:C_k_rips}.

\end{proof}


\section{Vertex degree in the k-NN graph}\label{sec:knn}

 We now demonstrate the broader application of the framework developed in Section \ref{sec:affine_inv}.
Let $G_k(\pointset)$ denote the $k$-nearest neighbor (k-NN) graph generated by a finite $\pointset\subset\R^d$.  
In other words, $x_1,x_2\in \pointset$ are connected by an edge if either $x_2$ is one of the $k$ nearest neighbors of $x_1$, or $x_1$ is one of the $k$ nearest nearest neighbors of $x_2$. 
Note that while above, $k$ denotes the homological degree, in this section only it refers to the parameter for constructing the graph.

When $\pointset$ is random, the degree of every vertex in $G_k(\pointset)$ is a random variable (larger or equal to $k$). In this section we show that the distribution of the vertex degree is universal. Note that for this functional, we only require the density to be upper bounded and  piecewise continuous, rather than the more restrictive conditions on good densities needed for the persistence measure.  

\begin{thm}\label{thm:univ_knn}
Let $\cP_\nn\sim\poisp{\nn}{f}$, and let $V_{\ell,\nn}$ be the number of vertices of degree $\ell$ in  $G_k(\cP_\nn)$. There exists a $K>k$, such that for all $k\le \ell\le K$, we have
\[
    \lim_{\nn\to\infty}\frac1\nn\mean{V_{\ell,\nn}} = p_\ell^*,
\]
where 
\[
    p_m^* = \prob{\deg(0;\cH_1\cup\set{0}) = \ell},
\]
and $\cH_1$ is a homogeneous Poisson process in $\R^d$ with rate $1$, so that 
$p_\ell^*$ depends on $d$, $k$, and $\ell$, but is independent of $f$ (i.e., universal). Furthermore, $\sum_{\ell=k}^K p_\ell^* = 1$.

Similarly, let $\cX_n\sim\binmp{n}{f}$, and let $V_{\ell,n}$ be the number of vertices of degree $\ell$ in $G_k(\cX_n$). Then almost surely, for all $k\le \ell\le K$,
\[
    \lim_{n\to\infty}\frac1\nn{V_{\ell,n}} =\lim_{n\to\infty}\frac1\nn\mean{V_{\ell,n}} = p_\ell^*.
\]
\end{thm}

To prove Theorem \ref{thm:univ_knn}, a key observation is that the neighborhood of every vertex in $G_k(\pointset)$ can be determined locally. To this end, 
we define
\eqb\label{eqn:cone}
    \cone_* := \set{x\in \R^d : \frac{\iprod{x,\bf 1}}{|x|} \ge \sqrt{3/2}},
\eqe
where ${\bf 1}= (1,0,\ldots,0)$,
so that $\cone_*$ is a cone whose  apex is the origin and whose angle is $\pi/3$.
Note that for every $x\in \R^d$, we can cover $\R^d$ with a finite collection of translated and rotated copies of $\cone_*$, denoted $\cone_1(x),\ldots,\cone_N(x)$. 

\begin{lem}\label{lem:knn_cones}
Let $\pointset\subset\R^d$ be locally-finite, and $G_k(\pointset)$ be the k-NN graph. 
Fix $x\in \pointset$, and define
\eqb\label{eqn:knn_R}
R(x;\pointset) := 2 \inf\set{r>0: |\pointset\cap B_r(x)\cap \cone_i(x)|\ge k+1,\ 1\le i \le N}.
\eqe
Then
\[
    \deg(x;\pointset) = \deg(x; \pointset\cap B_{R(x;\pointset)}(x)).
\]
In other words, the degree of $x$ is completely determined by the points lying in $B_{R(x;\pointset)}(x)$, and is not affected by changes to the point configuration outside this ball.
Additionally, taking $K=kN$, we have
\[
k\le\deg(x;\pointset) \le K.
\]
\end{lem}

\begin{proof}
We will prove the statement for $x=\origin$, assuming $\origin\in \pointset$, the result then follows from the translation and scale invariance of the k-NN graph.

Let $\cone_i=\cone_i(\origin)$, $1\le i\le N$, as defined above.
Note that for every $x_1,x_2\in \cone_i$, we have $|x_1-x_2| \le \max(|x_1|,|x_2|)$, since the angle between $x_1$ and $x_2$ is less than $\pi/3$ \eqref{eqn:cone}.
Let $\pointset\subset \R^d$, with $\origin\in \pointset$, and suppose that there exist $x_1,\ldots,x_m\in \cone_i\cap \pointset$, such that $0< |x_1|\le \cdots  \le |x_m|$, which implies that $|x_m-x_j| \le |x_m|$ for all $j\ne m$. If $x_m$ is connected to $\origin$ in $G_k(\pointset)$, then either: (a) $x_m$ is a $k$-nearest neighbor of $\origin$, which implies that $x_1,\ldots,x_{m-1}$ are all nearest neighbors as well; (b) $\origin$ is a $k$-nearest neighbor of $x_m$, but since $|x_m-x_j| \le |x_m|$, we conclude that $x_1,\ldots,x_{m-1}$ are all $k$-nearest neighbors of $x_m$ as well. In both cases we must have $m\le k$. In other words, each of the cones $\cone_i$ contains at most $k$ neighbors of $\origin$ in $G_k(\pointset)$.

Suppose that $R(\origin;\pointset)=r$. The definition \eqref{eqn:knn_R} implies that if $x$ is a neighbor of $\origin$ in $G_k(\pointset)$ then $x\in\cone_i\cap B_{r/2}(\origin)$ for some $i$. Additionally, note that $\cone_i$ contains at least $k$ points that are within distance $r/2$ from $x$. Therefore all the $k$-nearest neighbors of $x$ must lie in $B_{r/2}(\origin)$. This implies that $\deg(\origin;\pointset) = \deg(\origin;\pointset\cap B_{r}(\origin))$.

Finally, since each $\cone_i$ contains at most $k$ neighbors of $\origin$, we have $\deg(\origin;\pointset)\le kN$.
This concludes the proof.
\end{proof}

Using Lemma \ref{lem:knn_cones}, together with Theorems \ref{thm:univ_poiss} and \ref{thm:univ_binom}, we can now prove Theorem \ref{thm:univ_knn}.

\begin{proof}[Proof of Theorem \ref{thm:univ_knn}]
The functional we consider here is $\cH_\nn = V_{\ell,\nn}$, i.e., the number of vertices in $\cP_\nn$ with degree $\ell$. 
By definition, the k-NN graph is translation and scaling invariant, so \eqref{eqn:assum_translate} and \eqref{eqn:assum_scale} hold. In addition, the functional $\cH_\nn$ in this case has no explicit dependency in $\nn$ (only via $\cP_\nn$), and therefore \eqref{eqn:assum_homog} holds as well. For the binomial proof, we note that $V_{\ell,n}\le n$, so \eqref{eqn:assum_poly_1} holds.
Finally, since the degree of every vertex is between $k$ and $K$ (Lemma \ref{lem:knn_cones}), we have almost surely,
    \[
        |V_{\ell,n} - V_{\ell,n\pm\Delta}| \le (1+K)\Delta,
    \]
so that \eqref{eqn:assum_A_n}-\eqref{eqn:assum_poly_2} hold as well. Therefore, we are left with proving \eqref{eqn:assum_linear}, \eqref{eqn:assum_additive}, and \eqref{eqn:assum_stable}.

\noindent{\bf Proving linear scale \eqref{eqn:assum_linear}:}\\
Recall that here $\cP_\nn = \cP_\nn^*\sim\poisp{\nn}{\ind_{Q^d}}$. Using Mecke's formula, we have
\eqb\label{eqn:mean_V}
     \mean{V_{\ell,\nn}} = \nn\int_{Q^d}\prob{\deg(x;\cP_\nn^x) = \ell}dx,
\eqe
where $\cP_\nn^x := \cP^*_\nn\cup\set{x}$.
    
Fix $r\in (0,1/2)$, and $x\in [r,1-r]^d$. Recall the definition of $R(x;\cP_\nn^x)$ \eqref{eqn:knn_R}, and note that $R(x;\cP_\nn^x) \le r$ if and only if $R(x;\cP_\nn^x\cap B_{r}(x))\le r$. Also note that $\cP_\nn^* \cap B_{r}(x)$ is a homogeneous Poisson process with rate $\nn$. Therefore, using the scale and shift invariance of the k-NN graph, and taking  $\nn r^d=\lambda$, we have
\[
    \prob{R(x;\cP_\nn^x)>r} =  \probx{R(\origin;\cH_1^0)>\lambda^{1/d}},
\]
where $\cH_1$ is a homogeneous Poisson process in $\R^d$ with rate $1$. Similarly, from Lemma \ref{lem:knn_cones},
\[
\splitb
    \prob{\deg(x;\cP_\nn^x) = \ell, R(x;\cP_\nn^x) \le r} &= \prob{\deg(x;\cP_\nn^x\cap B_{r}(x)) = \ell, R(x;\cP_\nn^x) \le r}\\
    &= \prob{\deg(\origin;\cH_1^0\cap B_{\lambda^{1/d}}(\origin)) = \ell, R(\origin;\cH_1^0) \le \lambda^{1/d}} \\
    &= \prob{\deg(\origin;\cH_1^0)=\ell,R(\origin;\cH_1^0)\le \lambda^{1/d}}\\
    &= \prob{\deg(\origin;\cH_1^0)=\ell}- \prob{\deg(\origin;\cH_1^0)=\ell,R(\origin;\cH_1^0)> \lambda^{1/d}}.
\splite
\]
Denoting 
\[
p^*_\ell:=\prob{\deg(\origin;\cH_1^0)=\ell},
\]
we therefore have,
\[
p_\ell^* - \probx{R(\origin;\cH_1^0) > \lambda^{1/d}} \le \prob{\deg(x;\cP_\nn^x) = \ell} \le p_\ell^* + \probx{R(\origin;\cH_1^0)>\lambda^{1/d}}.
\]
Since $R(\origin;\cH_1^0)$ is almost surely finite (each  $\cone_i$ contains infinitely many points from $\cH_1$), we can take $\lambda$ large enough so that 
\[
\probx{R(\origin;\cH_1^0)>\lambda^{1/d}}< \eps.
\] 
Going back to \eqref{eqn:mean_V}, for every $x\in (0,1)^d$ and every $\eps>0$, we can take $\lambda$ and $\nn$ large enough so that $x\in [r,1-r]^d$, $\nn r^d = \lambda$, and
\[
|\prob{\deg(x;\cP_\nn^x) = \ell}-p_\ell^*| \le \eps,
\]
implying that
\[
\lim_{\nn\to\infty} \prob{\deg(x;\cP_\nn^x)=\ell} = p_\ell^*.
\]
Using the DCT for \eqref{eqn:mean_V} completes the proof.

\vspace{5pt}
\noindent{\bf Proving additivity \eqref{eqn:assum_additive}:}\\
Take any $M=m^d$, and define $\cP_\nn = \bigcup_i \cP_{\nn,i}$, where $\cP_{\nn,i}\sim\poisp{c_i\nn}{M\ind_{Q_i}}$ are independent. Denote by $V_{\ell,\nn}^{(i)}$ the number of degree-$\ell$ vertices generated by $\cP_{\nn,i}$. Our goal is to show that
\[
\lim_{\nn\to\infty} \frac1\nn \abs{\mean{V_{\ell,\nn}} - \sum_{i=1}^M \meanx{V_{\ell,\nn}^{(i)}}} = 0.
\]
Fix $r>0$, and let
\[
\cP_{\nn,i}^{(r)} = \set{x\in \cP_{\nn,i} : \inf_{z\in\partial Q_i}|x-z| \ge r}.
\]
Let $V_{\ell,\nn}^{(i,r)}$ be the number of degree-$\ell$ points $p$ in $G_k(\cP_{\nn,i})$, such that $p\in \cP_{\nn,i}^{(r)}$ and $R(p,\cP_{\nn,i}) \le r$. Then,
\[
V_{\ell,\nn}^{(i,r)} \le V_{\ell,\nn}^{(i)} \le V_{\ell,\nn}^{(i,r)} + \Delta_{1,i} + \Delta_{2,i},
\]
where $\Delta_{1,i}$ is the number of points $p\in \cP_{\nn,i}^{(r)}$ with $R(p;\cP_{\nn,i})>r$, and $\Delta_{2,i}$ is the number of points in $\cP_{\nn,i}\bs \cP_{\nn,i}^{(r)}$. 
Note that if $p\in\cP_{\nn,i}^{(r)}$ is such that $R(p,\cP_{\nn,i}) \le r$, then $\deg(p;\cP_{\nn,i}) = \deg(p;\cP_\nn)$. Therefore, we also have,
\[
\sum_{i=1}^M V_{\ell,\nn}^{(i,r)} \le V_{\ell,\nn} \le \sum_{i=1}^M V_{\ell,\nn}^{(i,r)} + \sum_{i=1}^M(\Delta_{1,i} + \Delta_{2,i}),
\]
which implies that
\[
\abs{\mean{V_{\ell,\nn}} - \sum_{i=1}^M \meanx{V_{\ell,\nn}^{(i)}}} \le \sum_{i=1}^M(\Delta_{1,i}+\Delta_{2,i}).
\]
Similarly to the steps in proving \eqref{eqn:assum_linear}, we can show that for all $\eps$, we can take $\lambda$ large enough, so that for  $\nn r^d = \lambda$,
\[
\limsup_{\nn\to\infty} \sum_{i=1}^M(\Delta_{1,i}+\Delta_{2,i}) < \eps,
\]
which proves \eqref{eqn:assum_additive}.

\noindent{\bf Proving continuity \eqref{eqn:assum_stable}:}\\
Let $f_1,f_2,\cdots$, be a sequence of  density functions as in \eqref{eqn:assum_stable}, and take $\cP_\nn = \cP_{\nn,i}\cup \Delta_{\nn,i}$. Denoting $V_{\ell,\nn}^{(i)}$ the number of degree-$\ell$ points in $G_k(\cP_{\nn,i})$, note that 
\[
    |V_{\ell,\nn} - V_{\ell,\nn}^{(i)}| \le (1+K)|\Delta_{\nn,i}|,
\]
since every point added from $\Delta_{\nn,i}$ can be of degree $\ell$, and can also  affect the degree of at most $K$ existing vertices in the graph.

Since $\mean{|\Delta_{\nn,i}|} = \delta_i\nn$, for every $\eps>0$ we can take $i$ large enough enough so that $\delta_i\le\eps/(1+K)$ which implies that
\[
\frac1\nn\meanx{|V_{\ell,\nn} - V_{\ell,\nn}^{(i)}|} \le \eps.
\]
This concludes the proof.

\end{proof}

\appendix

\section*{Appendix}

\section{Blaschke-Petkhanchin formula}

Our bounds for the number of critical faces for the \v Cech complex ($F_{k,\nn}(\lambda)$)  
 involve integrating over sets of points in $\R^d$ of size $k+1$. The  Blaschke-Petkahnchin formula provides a change of variables, extending the idea of polar coordinates, that significantly simplifies the integration in this context.
 
Let $\bx=(x_1,\ldots,x_{k+1}) \in (\R^d)^{k+1}$ be the coordinates we want to integrate over. Note that every $(k+1)$ points in general position, can be placed on a unique $(k-1)$-dimensional sphere in $\R^d$. Denoting this sphere by $\S(\bx)$, we then define the following
\[
\splitb
c(\bx) &:= \text{the center of } \S(\bx),\\
\rho(\bx) &:= \text{the radius of } \S(\bx),\\
\gamma(\bx) &:= \text{the $k$-dimensional plane containing }  \S(\bx),\\
\bth(\bx) &:= \text{the spherical coordinates of $\bx$ in } \S(\bx).
\splite
\]
The transformation $\bx\to (c(\bx), \rho(\bx), \gamma(\bx),\bth(\bx))$ is then a bijection, and we write $\bx = c+\rho(\gamma\circ\bth)$, 
where $\gamma\circ\bth$ takes the spherical coordinates $\bth$ on the plane $\gamma$.

Using this representation,
the following change of variables is known as the \emph{Blaschke-Petkantschin} (BP) formula. The version presented here was proved in \cite{edelsbrunner_expected_2017}.

\eqb\label{eqn:bp}
	\int_{(\R^d)^{k+1}} \varphi(\bx) d\bx = \int_{\R^d}\int_0^\infty \int_{\Gamma_{d,k}}\int_{(\S^{k-1})^{k+1}} \rho^{dk-1} \varphi(c+\rho(\gamma\circ \bth))(k! \vsimp(\bth))^{d-k+1}d\bth d\gamma d\rho dc,
\eqe
where $\Gamma(d,k)$ is the $k$-dimensional Grassmannian in $\R^d$.

\section{Strong law of large numbers}

The following lemma is inspired by the proof of  Theorem 3.17 in \cite{penrose_random_2003}, aiming to provide a more general context.

\begin{lem}\label{lem:exp_bounds}
    Let $f\in\cD(\cH)$, $\cX_n\sim\binmp{n}{f}$. Suppose that $\cH_n$ satisfies the conditions \eqref{eqn:assum_poly_1}-\eqref{eqn:assum_poly_2}. Then for every $\eps>0$ there exists $\gamma>0$, such that for $n$ large enough, we have
    \[
        \prob{\abs{\cH_n(\cX_n)-\mean{\cH_n(\cX_n)}} \ge \eps n} \le  e^{-n^\gamma}.
    \]
\end{lem}

\begin{proof}
  Define the martingale 
\[
M_i^{(n)} := \cmean{\cH_n(\cX_n)}{\cF_i},
\]
where $\cF_i = \sigma(\cX_i)$ (the $\sigma$-algebra generated by $\cX_i$), and the martingale difference
\[
D_i^{(n)} = M_i^{(n)} - M_{i-1}^{(n)}, \quad i\ge 1,
\]
 so that 
\[
\sum_{i=1}^n D_i^{(n)} = \cH_n(\cX_n) - \mean{\cH_n(\cX_n)} .
\]
Define $\cX_n^i = \cX_{n+1}\bs\set{X_i}$, and note that we can write
\[
D_i^{(n)} = \cmean{\cH_n(\cX_n) - \cH_n(\cX_n^i)}{\cF_i}.
\]
Recall the definition of $A_n$ in \eqref{eqn:assum_A_n}-\eqref{eqn:assum_poly_2}. We can define a similar event for $\cX_n^i$, denoted $\tilde A_{n,i}$, and take $A_{n,i} = A_n\cap \tilde A_{n,i}$. Note that from \eqref{eqn:assum_poly_2},  under the event $A_{n,i}$,
\[
\abs{\cH_n(\cX_n)-\cH_n(\cX_n^i)} \le \abs{\cH_n(\cX_n)-\cH_n(\cX_{n-1}^i)} + \abs{\cH_n(\cX^i_n)-\cH_n(\cX_{n-1}^i)} \le 2C_4n^b.
\]
Therefore, using \eqref{eqn:assum_poly_1},
\eqb\label{eqn:abs_D_i}
\splitb
    |D_i^{(n)}| &\le \cmean{\abs{\cH_n(\cX_n)-\cH_n(\cX_n^i)}\ind_{A_{n,i}}}{\cF_i} + \cmeanx{\abs{\cH_n(\cX_n)-\cH_n(\cX_n^i)}\ind_{A_{n,i}^c}}{\cF_i}\\
    &\le 2C_4n^b + C_3n^q\cprob{A_{n,i}^c}{\cF_i}.
\splite
\eqe

Denote $Z_{n,i} = \cprobx{A_{n,i}^c}{\cF_i}$, then
using \eqref{eqn:assum_A_n}, Markov's inequality, and the fact that $\probx{\tilde A_{n,i}} = \prob{A_n}$, yields
\[
\probx{Z_{n,i}> n^{-q}} \le n^q\prob{A_{n,i}^c} \le 2 n^q e^{-n^a}.
\]
From \cite{chalker_size_1999} (see also Theorem 2.9 in \cite{penrose_random_2003}), for any $\delta_1,\delta_2>0$, we have
\[
\prob{\abs{\sum_{i=1}^n D_i^{(n)}} > \delta_1} \le 2e^{-\frac{\delta_1^2}{32n\delta_2^2}} + \param{1+\frac{2\sup_i\|D_i^{(n)}\|_\infty}{\delta_1}} \sum_{i=1}^n \prob{|D_i^{(n)}
|>\delta_2}.
\]
Take $\delta_1 = \eps n$, and $\delta_2 = C_5n^b$ for some $C_5>2C_4$. Noting that $\|D_i^{(n)}\|_\infty \le C_3 n^q$, and that if $Z_{n,i}\le n^{-q}$ then from \eqref{eqn:abs_D_i} we have $|D_i^{(n)}| \le C_5 n^b$, we have for $n$ large,
\[
\splitb
\prob{|\cH_n(\cX_n)-\mean{\cH_n(\cX_n)}| \ge \eps n} &\le 2e^{-\frac{\eps^2n^2}{32 C_5^2 n^{2b+1}}} + (1+ 2\eps^{-1}C_3n^{q-1})\sum_{i=1}^n \probx{Z_{n,i}>n^{-q}},\\
&\le 2e^{-\frac{\eps^2}{32C_5^2 n^{2b-1}}} + 2(1+2\eps^{-1}C_3 n^{q-1})n^{q+1}e^{-n^a}\\
&\le e^{-n^{\gamma_1}} + e^{-n^{\gamma_2}},
\splite 
\]
where $\gamma_1 = (1-2b)/2>0$, and, $\gamma_2 = a/2$. Take any $\gamma < \min(\gamma_1,\gamma_2)$, then for $n$ large enough,
\[
\prob{|\cH_n(\cX_n)-\mean{\cH_n(\cX_n)}| \ge \eps n} \le e^{-n^\gamma},
\]
completing the proof.
\end{proof}

\section{Necessary conditions for criticality for Vietoris-Rips}
Here we prove a necessary condition for an edge in the Vietoris-Rips filtration to be critical, in the sense that the insertion of the edge and its star corresponds to some births or deaths in the persistence diagram (for any $k$).
Consider an edge $e=(x_1,x_2)$, and let $c = (x_1+x_2)/2$, and $\rho = |x_1-x_2|$. Recall that
\[
    I(x_1,x_2) := B_\rho(x_1)\cap B_\rho(x_2).
\]
Note, that
\[
   \sup_{I(x_1,x_2)} |x-c| = \frac{\sqrt{3}}{2}\rho.
\]
Therefore, if we take
\[
    \rho_{_\varobslash} := (1-\sqrt{3}/2)\rho,
\]
and
\[
I_{\varobslash}(x_1,x_2) := B_{\rho_{_\varobslash}}(c),
\]
then the distance between any $y\in I_{\varobslash}(x_1,x_2)$ and $z\in I(x_1,x_2)$ is less than $\rho$.
We use this observation  to prove the following lemma.

\begin{lem}\label{lem:forbidden}
Let $\cX\subset\R^d$, and $x_1,x_2\in \cX$.
If $e= (x_1,x_2)$ is critical in the Vietoris-Rips filtration generated by $\cX$, then  $I_{\varobslash}(x_1,x_2)\cap \cX=\emptyset$.
\end{lem}

\begin{proof}
Recall from Lemma \ref{lem:crit_rips}, that
\[
\splitb
F_k(e) &= \beta_{k-2}(\lk(e;\cX)),\quad k>1 \\
F_1(e) &=  \indf{\lk(e; \cX) = \emptyset},
\splite
\]
where the link is a Vietoris-Rips complex,
\[
\lk(e;\cX) = \cR_{\rho}(\cX'),\quad \cX' = \cX\cap I(x_1,x_2)\bs\{x_1,x_2\},\  \rho = |x_1-x_2|.
\]
Note that if $|\cX'| = 1$, then $F_k(e) = 0$ for all $k\ge 1$. Therefore, we assume $|\cX'| \ge 2$.

Clearly if $I_{\varobslash}(x_1,x_2)\cap \cX\ne\emptyset$ then $F_1(e) = 0$. 
To show that $F_k(e) = 0$ for $k>1$, suppose that  $y\in \cX\cap I_\varobslash(x_1,x_2)$. Then $|y-z|\le \rho$ for all $z\in \cX'$. 
Note that  for any simplex $\tau\in \lk(e;\cX)$ which does not include $y$, the simplex $\tau\cup\{y\}$ must also be in  $\lk(e;\cX)$, since all the vertices of $\tau$ are in $B_\rho(y)$. Hence, $\lk(e;\cX)$ is star-shaped with respect to $y$, and therefore 
contractible, and so 
$\beta_{k-2}(\lk(e;\cX)) = 0$ for all $k\ge 0$.

\end{proof}

\bibliography{zotero}
\bibliographystyle{plain}

\end{document}